# A Robust Traffic Control Model Considering Uncertainties in Turning Ratios

Hao Liu, Christian Claudel, Randy Machemehl, and Kenneth A. Perrine

***Abstract*— The effects of model parameter uncertainty on traffic flow control problems have recently drawn research attention. While the uncertainty in fundamental diagram related parameters has been investigated in the past, few articles have focused on network parameters uncertainty, including turning ratio uncertainty. To fill this gap, this article proposes a robust control model to deal with the uncertainties in the turning ratio by using distributionally robust chance constraints. The model allows one to compute the optimal control action that maximizes some objective, under all possible distributions of network parameters. We then apply this robust control framework to both a freeway network and an urban network, and evaluate the impact of uncertainty on optimal control inputs, over the test networks. The case studies show that compared to non-robust control, the proposed robust model can reduce congestion brought by the uncertainties and improve the overall throughput.**

## I. Introduction

TRAFFIC congestion has become a worldwide problem imposing a significant burden on both economy and environment. Traffic flow control is one of the primary methods to improve the efficiency of transportation systems, and a number of control methods have been proposed in the past decades for both freeway and urban networks.

The general goal of a traffic control method is to improve the average performance of the system, such as decreasing delay and increasing throughput. For most control methods, traffic flow model parameters, such as fundamental diagram, and external inputs, such as traffic densities, are assumed to be known and deterministic. However, most of these can be uncertain in practice, and neglecting these sources of uncertainty can lead to poor performance of the control scheme. The importance of considering randomness in traffic control problems has been widely recognized [1]–[4], and many efforts have been to handle the uncertainties.

For an urban network, it is commonly assumed all the vehicles travel at the same speed, which leads to parallel vehicle trajectories. As a result, the measure of effectiveness such as traffic delay is simplified. In reality, however, vehicle speed decreases when the traffic density reaches a certain degree of saturation. Fundamental diagrams are widely used to depict the speed change with traffic density. With traffic delay, uncertainties in vehicle arrival rates have a strong influence on control performance. A common way to study the effect of this is adding stochastic terms in the delay models. Heydecker [5] summarized the progress of this method. Based on queueing theory, Newell [6] developed a comprehensive approach to investigate the probabilistic arrivals. This method is based on steady-state analysis and may not be applicable for a short period study. To overcome this drawback, given the probability distribution functions (pdf) of arrival rates, Lo [7] proposed a phase clearance reliability (PCR) framework to investigate the probability of overflow during consecutive cycles, and this method was implemented on an adaptive signal control method for an arterial street [8].

For a freeway network, traffic flow is usually modeled by deterministic Partial Differential Equations (PDEs). The Lighthill-Whitham-Richards (LWR) [9], [10] model might be the most famous macroscopic traffic flow model, in which a fundamental diagram representing the relationship between traffic density and speed is required to obtain the solution. Besides the travel demand, the fundamental diagram parameters such as capacity and congested speed can be random since they are affected by external factors such as weather condition and driving behavior. The stochastic nature of the capacity has been studied broadly [11]–[13]. Furthermore, the initial density, as an input of optimization models, can be uncertain due to the sensor measurement errors. Based on the traffic flow control framework derived by Li *et al.* [14], [15], Liu *et al.* [16], [17] proposed a model to investigate the effect of the uncertainties in the initial densities on the control performance through chance constraints. Como *et al.* [18] proposed a traffic flow control method that is robust with respect to uncertainties in initial traffic volume and exogenous inflows.

As another concept to deal with the impact from uncertainties, resilience for dynamic traffic network has drawn considerable research attention [19]–[24]. Como *et al.* proposed distributed routing policies that are robust with respect to the reductions in traffic link capacities, and the robustness was evaluated in terms of both strong [20] and weak resilience [21]. Bianchin *et al.* [25] and Arnott *et al.* [26] demonstrated that networks adopting the suggestions from advanced data sharing techniques, such as Infrastructure-To-Vehicle (I2V) and Vehicle-To-Vehicle (V2V) communication technologies, do not guarantee the desired global performance. Yazıcıoğlu *et al.* [22], [23] proposed a resilient control for Hao Liu and Kenneth A. Perrine are with the Center for Transportation Research, The University of Texas at Austin, Austin, TX 78712 USA (e-mail: haoliu@utexas.edu; kperrine@utexas.edu).

Christian Claudel and Randy Machemehl are with the Department of Civil, Architectural and Environmental Engineering, The University of Texas at Austin, Austin, TX 78712 USA (e-mail: christian.claudel@utexas.edu; rbm@mail.utexas.edu).

dynamic networks by using variable speed limits to reduce systemic failures resulting from local routing decisions.

In addition to the sources mentioned in the literature above, the turning ratio is one of the main causes of uncertainty. However, to the authors' best knowledge, little effort has been put in developing robust control frameworks to handle such uncertainties to improve control reliability. To fill this gap, a second order cone program (SOCP) is proposed to study the impact of the uncertainties of random turning ratios on the traffic flow control. The proposed framework is based on the Lax-Hopf solution, derived by Mazaré *et al.* [27], to the LWR model. Unlike other analytical solutions such as front tracking method [28] which requires full knowledge of prior events and may have exponential growth of waves over time in some situations as waves "bounce" back and forth from the boundary conditions, this model is more efficient since it is grid-free and can calculate the traffic state directly from the initial and boundary conditions without any knowledge of prior events. In the proposed robust control model, we only assume that historical data is available, from which one can derive mean vectors and covariance matrices. These uncertainties are modeled as distributionally robust control chance constraints. Such constraints can be converted to SOC constraints [29] and solved by commercial SOCP solvers such as MOSEK [30].

The rest of this article is organized as follows. Section II reviews the derivation of the Lax-Hopf solution to the Hamilton-Jacobi (H-J) PDE, which is the building block of the proposed framework, and the constraints that the solutions need to satisfy. Following that, Section III shows the general form of deterministic traffic flow control models involving the constraints reviewed in Section II. Then, to handle the situation where turning ratios are random, a SOCP with distributionally robust chance constraints is proposed. The developed model is robust over all the distributions compatible with the moments (mean and covariance matrices) estimated from historical data. Sections IV and V implement the proposed model on both a freeway network and an urban network, and investigate the influence of the uncertainties on control inputs (on-ramp inflows and boundary inflows). Section VI summarizes the contribution and provides future research directions. Figure 1 serves as a reference to assist readers in grasping the role of each section and the structure of this article.

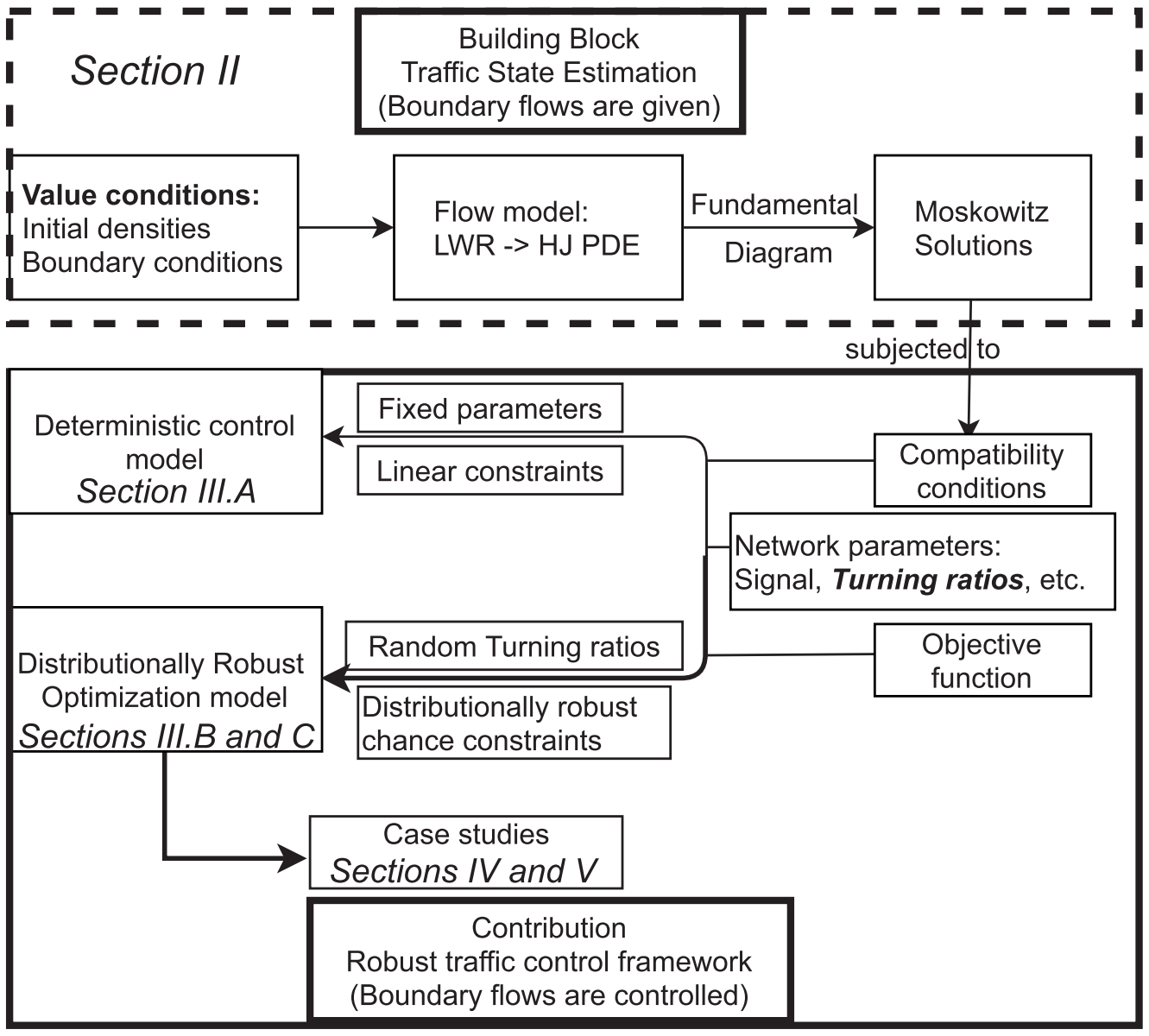

Fig. 1. Organization of the rest of this article.

## II. Review of the Lax-Hopf Solution

The proposed model is based on the Lax-Hopf solution [27], [31] of the H-J PDE. Li *et al.* [14], [15] has shown that a traffic flow model for a freeway network can be modeled as an optimization program with linear constraints based on the Lax-Hopf solution. As the building block of our robust control model, this section reviews the derivation of the solution and corresponding constraints. Part II-A introduces the traffic flow model, initial and boundary conditions; part II-B presents the Lax-Hopf solutions; part II-C shows the constraints that the Lax-Hopf solutions need to be satisfied with due to the compatibility conditions. For details regarding the derivation and proof, the readers are referred to [14], [32].

### A. H-J PDE

The LWR model [9], [10] is a widely used macroscopic model depicting relationships among traffic flow characteristics

$$\frac{\partial \rho(t,x)}{\partial t} + \frac{\partial \psi(\rho(t,x))}{\partial x} = 0, \tag{1}$$

where $\rho(t,x)$ is the density of the point $x$ away from a reference point at time $t$, $\psi$ denotes the experimental relationship, which is defined as the fundamental diagram between flow and density. The proposed model is applicable for any concave and piecewise linear fundamental diagram. For non-linear FDs, in order to employ the proposed model, we need to obtain their piecewise linear approximations. For simplicity, a triangular fundamental diagram is utilized in this article,

$$\psi(\rho) = \begin{cases} v_f \rho & \rho \in [0, \rho_c] \\ w(\rho - \rho_m) & \rho \in [\rho_c, \rho_m], \end{cases} \tag{2}$$

where $v_f$ is the free flow speed, $w$ is the congestion speed, $\rho_c$ is the critical density where the flow reaches its capacity, $\rho_m$ is the jam density, where the flow is zero due to the total congestion. Alternatively, by integrating the LWR PDE in space, another traffic flow model, H-J PDE, can be expressed as

$$\frac{\partial M(t,x)}{\partial t} - \psi(-\frac{\partial M(t,x)}{\partial x}) = 0, \tag{3}$$

where $M(t,x)$, known as the Moskowitz function [33], denotes the index of the vehicles at point $(t,x)$. To solve this function, the spatial domain $[\xi, \chi]$ is divided evenly into $k_{max}$ segments and the time domain $[0, t_{max}]$ is divided evenly into $n_{max}$ segments. Let $K = \{1, \ldots, k_{max}\}$ and $N = \{1, \ldots, n_{max}\}$. Assuming the initial density in each spatial segment and the flow in each time step are constant, the piecewise affine initial condition $M_k(t,x)$, upstream boundary condition

$\gamma_n(t,x)$, and downstream boundary condition $\beta_n(t,x)$ are defined as

$$M_k(t,x)=\begin{cases}-\sum_{i=1}^{k-1}\rho(i)X & \\ -\rho(k)(x-(k-1)X), & \text{if } t=0 \\ & \text{and } x\in[(k-1)X,kX] \\ +\infty, & \text{otherwise}\end{cases} \tag{4}$$

$$\gamma_n(t,x)=\begin{cases}\sum_{i=1}^{n-1}q_{\text{in}}(i)T & \\ +q_{\text{in}}(n)(t-(n-1)T), & \text{if } x=\xi \\ & \text{and } t\in[(n-1)T,nT] \\ +\infty, & \text{otherwise}\end{cases} \tag{5}$$

$$\beta_n(t,x)=\begin{cases}\sum_{i=1}^{n-1}q_{\text{out}}(i)T & \\ +q_{\text{out}}(n)(t-(n-1)T) & \\ -\sum_{k=1}^{k_{max}}\rho(k)X, & \text{if } x=\chi \\ & \text{and } t\in[(n-1)T,nT] \\ +\infty, & \text{otherwise}\end{cases} \tag{6}$$

where $X$ and $T$ are the spatial segment length and time step size, respectively, $\rho(i)$ is the initial density for the $i$th spatial segment and $q_{\text{in}}(i)$ and $q_{\text{out}}(i)$ are the inflow and outflow for the $i$th time step, respectively.

### B. Lax-Hopf Solutions

The Barron-Jensen/Frankowska (B-J/F) solution [34], [35] fully characterized by the Lax-Hopf formula was incorporated to solve the H-J equation.

*Definition 1 (Value Condition):* A value condition $c(\cdot,\cdot)$ is a lower semicontinuous function defined on a subset of $[0,t_{max}]\times[\xi,\chi]$.

The initial conditions and boundary conditions are regarded as value conditions.

*Proposition 1 (Lax-Hopf Formula):* Let $\psi(\cdot)$ be a concave and continuous Hamiltonian, and let $c(\cdot,\cdot)$ be a value condition. The B-J/F solution $M_c(\cdot,\cdot)$ to (3) associated with $c(\cdot,\cdot)$ is defined [36]–[38] by

$$M_c(t,x)=\inf_{(u,T)\in(\varphi^*)\times R_+}(c(t-T,x+Tu)+T\varphi^*(u)) \tag{7}$$

where $\varphi^*(\cdot)$ is the Legendre-Fenchel transform of an upper semicontinuous Hamiltonian $\psi(\cdot)$, which is given by,

$$\varphi^*(u):=\sup_{p\in Dom(\psi)}[p\cdot u+\psi(p)] \tag{8}$$

Based on this proposition, the Moskowitz solution from value conditions (4)-(6) can be expressed as (9)-(11). Many formulas are divided into multiple lines, and the commas indicate the end position of formulas. For readers interested in the derivation, see [27] and [31] for more on the derivation.

$$M_{M_k}(t,x)=\begin{cases}+\infty, & \text{if } x\le(k-1)X+tw \text{ or } x\ge kX+v_ft \quad (9a)\\ -\sum_{i=1}^{k-1}\rho(i)X+\rho(k)(tv_f+(k-1)X-x), & \text{if } x\ge(k-1)X+v_ft \text{ and } x\le kX+v_ft \text{ and } \rho(k)\le\rho_c \quad (9b)\\ -\sum_{i=1}^{k-1}\rho(i)X+\rho_c(tv_f+(k-1)X-x), & \text{if } x\le(k-1)X+v_ft \text{ and } x\ge(k-1)X+tw \text{ and } \rho(k)\le\rho_c \quad (9c)\\ -\sum_{i=1}^{k-1}\rho(i)X+\rho(k)(tw+(k-1)X-x)-\rho_m tw, & \text{if } x\le kX+tw \text{ and } x\ge(k-1)X+tw \text{ and } \rho(k)\ge\rho_c \quad (9d)\\ -\sum_{i=1}^{k}\rho(i)X+\rho_c(tw+kX-x)-\rho_m tw, & \text{if } x\le kX+tv_f \text{ and } x\ge kX+tw \text{ and } \rho(k)\ge\rho_c \quad (9e)\end{cases}$$

$$M_{\gamma_n}(t,x)=\begin{cases}+\infty, & \text{if } t\le(n-1)T+\frac{x-\xi}{v_f} \quad (10a)\\ \sum_{i=1}^{n-1}q_{\text{in}}(i)T+q_{\text{in}}(n)(t-\frac{x-\xi}{v_f}-(n-1)T), & \text{if } t\ge(n-1)T+\frac{x-\xi}{v_f} \text{ and } t\le nT+\frac{x-\xi}{v_f} \quad (10b)\\ \sum_{i=1}^{n}q_{\text{in}}(i)T+\rho_c v_f(t-\frac{x-\xi}{v_f}-nT), & \text{otherwise} \quad (10c)\end{cases}$$

$$M_{\beta_n}(t,x)=\begin{cases}+\infty, & \text{if } t\le(n-1)T+\frac{x-\chi}{w} \quad (11a)\\ -\sum_{k=1}^{k_{max}}\rho(k)X+\sum_{i=1}^{n-1}q_{\text{out}}(i)T+q_{\text{out}}(n)(t-\frac{x-\chi}{w}-(n-1)T)-\rho_m(x-\chi), & \text{if } t\ge(n-1)T+\frac{x-\chi}{w} \text{ and } t\le nT+\frac{x-\chi}{w} \quad (11b)\\ -\sum_{k=1}^{k_{max}}\rho(k)X+\sum_{i=1}^{n}q_{\text{out}}(i)T+\rho_c v_f(t-nT-\frac{x-\chi}{v_f}), & \text{otherwise} \quad (11c)\end{cases}$$

### C. Linear Constraints

The Moskowitz solutions (9)-(11) show that each value condition generates one solution at a certain point in the domain of value conditions. The corresponding compatibility conditions need to be satisfied by these solutions.

The Lax-Hopf formula (7) leads to the inf-morphism property [36].

*Proposition 2 (Inf-Morphism Property):* Let the value condition $\boldsymbol{c}(\cdot,\cdot)$ be minimum of a finite number of lower semi-continuous functions:

$$\forall(t,x) \in [0, t_{max}] \times [\xi, \chi], \quad \boldsymbol{c}(t,x) := \min_{j\in J} \boldsymbol{c}_j(t,x) \tag{12}$$

The corresponding solution $\boldsymbol{M_c}(\cdot,\cdot)$ can be decomposed [36], [37] as

$$\forall(t,x) \in [0, t_{max}] \times [\xi, \chi], \quad \boldsymbol{M_c}(t,x) := \min_{j\in J} \boldsymbol{M_{c_j}}(t,x) \tag{13}$$

Based on the *Inf-morphism* property, the Moskowitz solutions (9)-(11) have to satisfy the compatibility conditions [31].

*Proposition 3 (Compatibility Conditions):* Use the value condition $\boldsymbol{c}(t,x)$ and the corresponding solution in *Proposition 2*. The equality $\forall(t,x) \in Dom(\boldsymbol{c}), \boldsymbol{M_c}(t,x) = \boldsymbol{c}(t,x)$ is valid if and only if the inequalities below are satisfied,

$$\boldsymbol{M_{c_j}}(t,x) \geq \boldsymbol{c}_i(t,x), \quad \forall(t,x) \in Dom(\boldsymbol{c}_i), \ \forall(i,j) \in J^2 \tag{14}$$

In detail, these constraints can be expanded as [14], [32], (15)–(17), as shown at the bottom of the page. The Moskowitz solutions in these constraints are piecewise linear function of inflows and outflows. A traffic flow control model needs to satisfy these constraints to make the problem compatible. Unlike other traffic flow control methods [39]–[41] in which the PDEs are discretized to ODEs to employ available algorithms, such as gradient descent [42], to solve the optimization model, this framework does not require any discretization or approximation of the corresponding PDE.

## III. Robust Model With Distributionally Robust Chance Constraints

This section first shows the general form of deterministic traffic flow control models in which model parameters are constant. Then, the uncertainties in the turning ratios are raised and the robust control model is developed to deal with such randomness.

### A. Deterministic Control Model

Founded on the Lax-Hopf solution, a traffic flow control model for a freeway network can be expressed as,

$$\begin{aligned} &\min \ f(\boldsymbol{c},\boldsymbol{x}) \\ &s.t. \ (15)-(17), \quad \forall l \in L \\ &\begin{bmatrix} \boldsymbol{q}^z_{\text{out}} \\ q^z_{\text{off}} \end{bmatrix} = \begin{bmatrix} \boldsymbol{P}^z_{\boldsymbol{1}} & \boldsymbol{P}^z_{\boldsymbol{2}} \\ \boldsymbol{P}^z_{\boldsymbol{3}} & 0 \end{bmatrix} \begin{bmatrix} \boldsymbol{q}^z_{\text{in}} \\ q^z_{\text{on}} \end{bmatrix}, \quad \forall z \in Z \end{aligned} \tag{18}$$

where $l$ and $z$ are the index of links and nodes, and $L$ and $Z$ are the sets of links and nodes. The decision variable $\boldsymbol{x}$ vector is

$$\boldsymbol{x} := \{q_{\text{in}}(i,j), q_{\text{out}}(i,j) : i \in N, j \in L\}, \tag{19}$$

where $N$ is the set of time steps, $q_{\text{in}}(i,j)$ and $q_{\text{out}}(i,j)$ are the inflow and outflow of link $j$ at time $i$, respectively. $f(\boldsymbol{c},\boldsymbol{x})$

$$\begin{cases} M_{M_k}(0,x_p) \geq M_p(0,x_p) & \forall(k,p) \in K^2 \\ M_{M_k}(pT,\chi) \geq \beta_p(pT,\chi) & \forall k \in K, \quad \forall p \in N \\ M_{M_k}(\dfrac{\chi - x_k}{v_f},\chi) \geq \beta_p(\dfrac{\chi - x_k}{v_f},\chi) & \forall k \in K, \quad \forall p \in N \\ \qquad \text{s.t.} \ \dfrac{\chi - x_k}{v_f} & \in [(p-1)T, pT] \\ M_{M_k}(pT,\xi) \geq \gamma_p(pT,\xi) & \forall k \in K, \quad \forall p \in N \\ M_{M_k}(\dfrac{\xi - x_{k-1}}{w},\xi) \geq \gamma_p(\dfrac{\xi - x_{k-1}}{w},\xi) & \forall k \in K, \quad \forall p \in N \\ \qquad \text{s.t.} \ \dfrac{\xi - x_{k-1}}{w} & \in [(p-1)T, pT] \end{cases} \tag{15}$$

$$\begin{cases} M_{\gamma_n}(pT,\xi) \geq \gamma_p(pT,\xi) & \forall(n,p) \in N^2 \\ M_{\gamma_n}(pT,\chi) \geq \beta_p(pT,\chi) & \forall(n,p) \in N^2 \\ M_{\gamma_n}(nT + \dfrac{\chi - \xi}{v_f},\chi) \geq \beta_p(nT + \dfrac{\chi - \xi}{v_f},\chi) & \forall(n,p) \in N^2 \\ \qquad \text{s.t.} \ nT + \dfrac{\chi - \xi}{v_f} & \in [(p-1)T, pT] \end{cases} \tag{16}$$

$$\begin{cases} M_{\beta_n}(pT,\xi) \geq \gamma_p(pT,\xi) & \forall(n,p) \in N^2 \\ M_{\beta_n}(nT + \dfrac{\xi - \chi}{w},\xi) \geq \gamma_p(nT + \dfrac{\xi - \chi}{w},\xi) & \forall(n,p) \in N^2 \\ \qquad \text{s.t.} \ nT + \dfrac{\xi - \chi}{w} & \in [(p-1)T, pT] \\ M_{\beta_n}(pT,\chi) \geq \beta_p(pT,\chi) & \forall(n,p) \in N^2 \end{cases} \tag{17}$$

is a general form of the objective function, and $\boldsymbol{c}$ is the involved parameters. The first constraint means the inflows and outflows of all links need to satisfy the compatibility condition; the second constraint is the node model representing the flow transition. Note that a node has at least more than one incoming link (including on-ramps) or more than one outgoing link (including off-ramps).

Let $N_{out}^z$ and $N_{in}^z$ represent the number of outgoing links and incoming links at node $z$. $\boldsymbol{q}_{\text{in}}^z$ and $\boldsymbol{q}_{\text{out}}^z$ are two column vectors denoting the incoming flows and outgoing flows at node $z$, respectively; $q_z^{\text{on}}$ and $q_z^{\text{off}}$ are two scalars representing the on-ramp and off-ramp flows; $\boldsymbol{P}_{\mathbf{1}}^z$ is a $N_{out}^z \times N_{in}^z$ matrix of which each element $\boldsymbol{P}_{\mathbf{1}}^z(i, j)$ means the proportion of the vehicles from the $j$th incoming link going into the $i$th outgoing link at node $z$; $\boldsymbol{P}_{\mathbf{2}}^z$ is a column vector with dimension of $N_{out} \times 1$ of which each element $\boldsymbol{P}_{\mathbf{2}}^z(i)$ means the proportion of the vehicles from on-ramp going into the $i$th outgoing link; $\boldsymbol{P}_{\mathbf{3}}^z$ is a row vector with dimension of $1 \times N_{in}$ of which each element $\boldsymbol{P}_{\mathbf{3}}^z(j)$ means the proportion of the vehicles from the $j$th incoming link departing from the off-ramp. In addition, we assume no vehicles coming from an on-ramp would depart from the off-ramp at the same node, which makes the last element in the transition matrix equal to 0.

### B. Introduction of Distributionally Robust Chance Constraints

In reality, the turning ratio matrices are not always deterministic, and only prior information of their distributions such as moments can be extracted from historical data. Under this situation, the distributionally robust optimization model, in short DRO referred by Delage and Ye [43], is a proper method to study the effect of such uncertainties. This modeling framework has received considerable attention in research communities such as operations research and machine learning. A comprehensive review of DRO can be found in [44].

Let $\tilde{\boldsymbol{\xi}} \in R^k$ denote the random parameters, its ambiguity set is defined as the set of distributions that are consistent with the prior knowledge about the uncertainty. Assume $\boldsymbol{\mu}$ and $\boldsymbol{\Gamma}$ are its expectation and covariance matrices, and they are the only information known. Then, its ambiguity set can be expressed as

$$\mathcal{P} = \{\mathbb{P} \in P(R^k) : \mathbb{E}_{\mathbb{P}}[\tilde{\boldsymbol{\xi}}] = \boldsymbol{\mu}, \quad \mathbb{E}_{\mathbb{P}}[(\tilde{\boldsymbol{\xi}} - \boldsymbol{\mu})(\tilde{\boldsymbol{\xi}} - \boldsymbol{\mu})^T] = \boldsymbol{\Gamma}\}. \tag{20}$$

The goal of a DRO is to optimize the worst-case objective value over the ambiguity set. For example, the objective function of a stochastic program in which the random parameters' distributions are known can be expressed as

$$\min_{\boldsymbol{x}} \; R[f(\boldsymbol{x}, \tilde{\boldsymbol{\xi}})] \tag{21}$$

where $\boldsymbol{x}$ is the decision variable vector, $R$ is the risk measure, such as expectation and Value at Risk (VaR). When an ambiguity set of $\tilde{\boldsymbol{\xi}}$ is given, this model can be transformed as a DRO

$$\min_{\boldsymbol{x}} \; \max_{\mathbb{P} \in \mathcal{P}} \; R_{\mathbb{P}}[f(\boldsymbol{x}, \tilde{\boldsymbol{\xi}})]. \tag{22}$$

On the other hand, if the randomness is involved in constraints, a common way to develop robust model is to replace the deterministic constraints with chance constraints

$$P[h(\boldsymbol{x}, \tilde{\boldsymbol{\xi}}) \geq 0] \geq 1 - \alpha \tag{23}$$

which indicates that the constraint $h(\boldsymbol{x}, \tilde{\boldsymbol{\xi}}) \geq 0$ hold with confidence level of $1 - \alpha$. Similarly, if only the ambiguity set is known, the distributionally robust chance constraint can be expressed as

$$\min_{\mathbb{P} \in \mathcal{P}} P[h(\boldsymbol{x}, \tilde{\boldsymbol{\xi}}) \geq 0] \geq 1 - \alpha \tag{24}$$

which means the minimum of the probability, i.e. the worst case, that the constraint holds under all possible distributions is larger than the confidence level.

If $h(\boldsymbol{x}, \tilde{\boldsymbol{\xi}})$ is an affine function of $\boldsymbol{x}$, then the distributionally robust chance constrain can be expressed as

$$\min_{\mathbb{P} \in \mathcal{P}} P[\tilde{\boldsymbol{a}}^T \boldsymbol{x} + \tilde{b} \leq 0] \geq 1 - \alpha. \tag{25}$$

Let $\boldsymbol{d} = [\tilde{\boldsymbol{a}}^T, \tilde{b}]^T$. If $\mathcal{P}$ is its ambiguity set with known expectation $\hat{\boldsymbol{d}}$ and covariance matrix $\boldsymbol{\Gamma}$ and $1 - \alpha > 0.5$, (25) can be converted to a convex second-order cone (SOC) constraint [29]

$$\kappa_\alpha \sigma(\tilde{x}) + \hat{\varphi}(\tilde{x}) \leq 0, \quad \kappa_\alpha = \sqrt{(1 - \alpha)/\alpha}, \tag{26}$$

where $\tilde{\boldsymbol{x}} = [\boldsymbol{x}^T, 1]^T$, $\hat{\varphi}(\tilde{\boldsymbol{x}}) = \hat{\boldsymbol{d}}^T \tilde{\boldsymbol{x}}$, $\boldsymbol{\sigma}^2(\tilde{\boldsymbol{x}}) = \tilde{\boldsymbol{x}}^T \boldsymbol{\Gamma} \tilde{\boldsymbol{x}}$ and $\boldsymbol{\sigma}(\tilde{\boldsymbol{x}}) = \|\boldsymbol{\Gamma}^{\frac{1}{2}} \tilde{\boldsymbol{x}}\|_2$. The proof can be found in [29].

### C. Robust Control Model as a SOCP

The turning ratios are involved in the second constraint in (18). However, we cannot convert these equality constraints to chance constraints directly due to the fact that for any feasible solution, the probability an equality constraint holds is always zero if the distribution is continuous. Therefore, we need to transform those constraints to inequality form before adding chance constraints. To this end, we make following definitions:

1. The on-ramps and off-ramps are regarded as incoming links and outgoing links, respectively;
2. The links are divided into two groups: incoming boundary links and other links. Incoming boundary links are the links through which the vehicles flow into the network. For example, links 1 and 4 and all on-ramps in Figure 2 are incoming boundary links. Let $L$ denote the set of links and $L_{\text{in}}$ be the set of incoming boundary links;
3. All the nodes contain more than one incoming link or more than one outgoing link. Otherwise, the turning ratio is always 1. Let $Z$ be the set of nodes, and $L_{\text{in}}^z$ and $L_{\text{out}}^z$ be the incoming link set and outgoing link set of node $z$;
4. Let $\boldsymbol{P}^z$ be a $m_z \times n_z$ matrix where $m_z$ and $n_z$ are the number of outgoing links and incoming links at node $z$. $P^z(i, j)$ represents the ratio of vehicles from link $L_{\text{in}}^z(j)$ to link $L_{\text{out}}^z(i)$;
5. The inflows except for the incoming boundary links can be expressed as: $q_{\text{in}}(i, j) = \sum_{k \in L_{\text{in}}^z} P^z(j, k) q_{\text{out}}(i, k)$, $\forall j \in L_{\text{out}}^z$.

By this way, we can replace most of the $q_{\text{in}}(i, j)$'s with a function of $q_{\text{out}}(i, j)$'s, and the equality constraint in (18) can be removed. The new decision variable $\boldsymbol{x}$ is defined as follows:

$$\boldsymbol{x} := \{q_{\text{in}}(i, j) : i \in N, j \in L_{\text{in}}\} \cup \{q_{\text{out}}(i, j) : i \in N, j \in L\}. \tag{27}$$

In addition, for the upstream boundary conditions (5) and Moskowitz solutions (9) related to $q_{\text{in}}(i, j)$'s, we need to reformat them as, (28) and (29), as shown at the bottom of the page, where $l$ is the link index, $z$ is the node from which link $l$ starts, i.e. link $l$ is one of the outgoing links of node $z$. All other initial conditions and Moskowitz solutions will not be changed. Now, the $\gamma_n$'s and $M_{\gamma_n}$'s in the constraints (15)-(17) are replaced by (28) and (29). The new constraints are linear functions of $\boldsymbol{q}_{\text{out}}$ with random coefficient $\boldsymbol{P}$. Therefore, the chance constraint and distributionally robust chance constraint of each of them can be converted to a SOC constraint by following the conversion from (25) to (26).

In the rest of this section, for simplicity, we only derive the SOC conversion constraint the fourth constraint in (15):

$$\min_{\mathbb{P}\in\mathcal{P}} P[M_{M_k}(pT, \xi) \geq \gamma_p(pT, \xi)] \geq 1-\alpha \quad \forall k \in K, \forall p \in N. \tag{30}$$

The remaining SOC constraints for (15)-(17) can be obtained in the same way, and the expressions are shown as (40), (47) and (57) in the appendix. Each SOC constraint forces the optimal solutions to be feasible for $100(1-\alpha)$ percent of time while the optimal inflows or outflows may exceed their upper limits for $100\alpha$ percent of time. In the latter case, the excess portion of the inflows cannot enter the link at the expected time step and may lead to congestion. In other words, each SOC constraint ensures the LWR PDE (1) holds true at the domain of value conditions for the optimal solutions with probability $1-\alpha$.

Let $\tilde{\boldsymbol{P}}^z$ be the random turning ratio matrix at node $z$, let $\boldsymbol{P}^z$ and $\boldsymbol{\Gamma}$ be its the expectation and covariance matrix. First, we need to replace the inflow of an outgoing link at node $z$, $L^z_{\text{out}}(j)$, with a function of the outflows of incoming links at node $z$, shown as

$$q_{\text{in}}(i, L^z_{\text{out}}(j)) = \sum_{r\in L^z_{\text{in}}} \tilde{P}^z(L^z_{\text{out}}(j), r) q_{\text{out}}(i, r), \tag{31}$$

where $q_{\text{in}}(i, L^z_{\text{out}}(j))$ is the inflow of link $L^z_{\text{out}}(j)$ at time $i$, $q_{\text{out}}(i, r)$ is the outflow of link $r$ at time $i$, $\tilde{P}^z(L^z_{\text{out}}(j), r)$ is the ratio of vehicles from link $r$ to link $L^z_{\text{out}}(j)$, and $L^z_{\text{in}}$ is the set of incoming links at node $z$. By substituting (31) into (30), we obtain

$$\sum_{i=1}^{p}\Big(\sum_{r\in L^z_{\text{in}}} q_{\text{out}}(i, r)\tilde{P}^z(L^z_{\text{out}}(j), r)\Big)T - M_{M_k}(pT, \xi) \leq 0, \quad \forall k \in K, \quad \forall p \in N. \tag{32}$$

In this constraint, the decision variables are $q_{\text{out}}(i, r)$'s, and $i = 1, 2, \ldots, p$ and $r = L^z_{\text{in}}(1), L^z_{\text{in}}(2), \ldots, L^z_{\text{in}}(|\boldsymbol{L}^z_{\text{in}}|)$. Note that $M_{M_k}(.)$ is the Moskowitz solution from initial densities and dependent on turning ratios, so we do not need to reformat this part. Therefore, based on (26), the decision variable vector can be expressed as

$$\begin{aligned}\tilde{\boldsymbol{x}}^p_1 := [&\overbrace{q_{\text{out}}(1, L^z_{\text{in}}(1)), q_{\text{out}}(2, L^z_{\text{in}}(1)), \ldots, q_{\text{out}}(p, L^z_{\text{in}}(1))}^{p},\\ &\overbrace{\cdots\cdots}^{(n_z-2)\times p},\\ &\overbrace{q_{\text{out}}(1, L^z_{\text{in}}(n_z)), q_{\text{out}}(2, L^z_{\text{in}}(n_z)), \ldots, q_{\text{out}}(p, L^z_{\text{in}}(n_z))}^{p},\\ &1]^T,\end{aligned} \tag{33}$$

$$\gamma_n(t, x, l) = \begin{cases} \sum_{i=1}^{n-1}\big(\sum_{k\in L^z_{\text{in}}} q_{\text{out}}(i, k)P^z(l, k)\big)T \\ +\big(\sum_{k\in L^z_{\text{in}}} q_{\text{out}}(n, k)P^z(l, k)\big)(t-(n-1)T), & \text{if } x = \xi \text{ and } t \in [(n-1)T, nT] \\ +\infty, & \text{otherwise} \end{cases} \tag{28}$$

and

$$M_{\gamma_n}(t, x, l) = \begin{cases} +\infty, & \text{if } t \leq (n-1)T + \frac{x-\xi}{v_f} \\ \sum_{i=1}^{n-1}\big(\sum_{k\in L^z_{\text{in}}} q_{\text{out}}(i, k)P^z(l, k)\big)T + \\ \big(\sum_{k\in L^z_{\text{in}}} q_{\text{out}}(n, k)P^z(l, k)\big) & \text{if } t \geq (n-1)T + \frac{x-\xi}{v_f} \\ (t - \frac{x-\xi}{v_f} - (n-1)T), & \text{and } t \leq nT + \frac{x-\xi}{v_f} \\ \sum_{i=1}^{n}\big(\sum_{k\in L^z_{\text{in}}} q_{\text{out}}(i, k)P^z(l, k)\big)T + \rho_c v_f & \text{otherwise} \\ (t - \frac{x-\xi}{v_f} - nT) \end{cases} \tag{29}$$

Then, $\hat{\boldsymbol{d}}_1^p$ is the coefficient of each decision variable in (32), and it can be expressed as

$$\begin{aligned}\hat{\boldsymbol{d}}_1^p := &[\overbrace{T\boldsymbol{P}^z(L^z_{\text{out}}(i), L^z_{\text{in}}(1)), \ldots, T\boldsymbol{P}^z(L^z_{\text{out}}(i), L^z_{\text{in}}(1)),}^{p} \\ &\overbrace{\cdots\cdots}^{(n_z-2)\times p}, \\ &\overbrace{T\boldsymbol{P}^z(L^z_{\text{out}}(i), L^z_{\text{in}}(n_z)), \ldots, T\boldsymbol{P}^z(L^z_{\text{out}}(i), L^z_{\text{in}}(n_z)),}^{p} \\ &-M_{M_k}(pT, \xi)], \end{aligned} \tag{34}$$

where each block of $\tilde{\boldsymbol{x}}_1^p$ and $\hat{\boldsymbol{d}}_1^p$ indicated by the overbrace is the outflows for the same incoming link at $p$ time steps and the corresponding coefficients. Based on the dimension of the decision variable vector, the covariance matrix is a $p \times (n_z+1)$ matrix. We assume the uncertainties of the turning ratios do not change over time steps, so

$$\text{cov}(\hat{\boldsymbol{d}}_1^p(i), \hat{\boldsymbol{d}}_1^p(j)) = \begin{cases} \Gamma(n,n) & \text{if } i,j \in \text{block n} \\ \Gamma(n,m) & \text{if } i \in \text{block n}, j \in \text{block m} \\ 0 & \text{if } i = n_z+1 \text{ or } j = n_z+1 \end{cases} \tag{35}$$

Therefore, the integral covariance matrix (36), as shown at the bottom of the page.

Therefore, by substituting Equations (33), (34) and (36) to Equation (26), we can convert the optimization model (18) to a program with SOC constraints. If the objective function is linear, the model becomes a SOCP.

## IV. Case Study on a Freeway Network

This section implements the proposed framework on a network to test the impact of uncertainties in turning ratios on the control results.

### A. Network and Problem Settings

The network employed is shown in Figure 2. This network consists of 6 links, 3 nodes (excluding the boundary nodes), 2 on-ramps and 2 off-ramps. Links 1 and 4 are incoming links. Each link has a length of 1.2km and is divided into 2 segments; the simulation time is 500s which is made into 25 equal time steps. Consider the free flow speed $v_f$ = 30m/s, critical density $\rho_c = 0.0175$ veh/m/lane, capacity $C = 0.5250$ veh/s/lane, jam density $\rho_m = 0.2250$ veh/m/lane, and the congestion speed $w = -5.5$m/s. Links 1-3 have 4 lanes, and links 4-6 have 3 lanes. The related sets, vectors and matrices are defined as follows:

1. Link sets: $L = \{i \in \mathbb{Z} : 1 \le i \le 6\}$, $L_{\text{in}} = \{1, 4\}$;

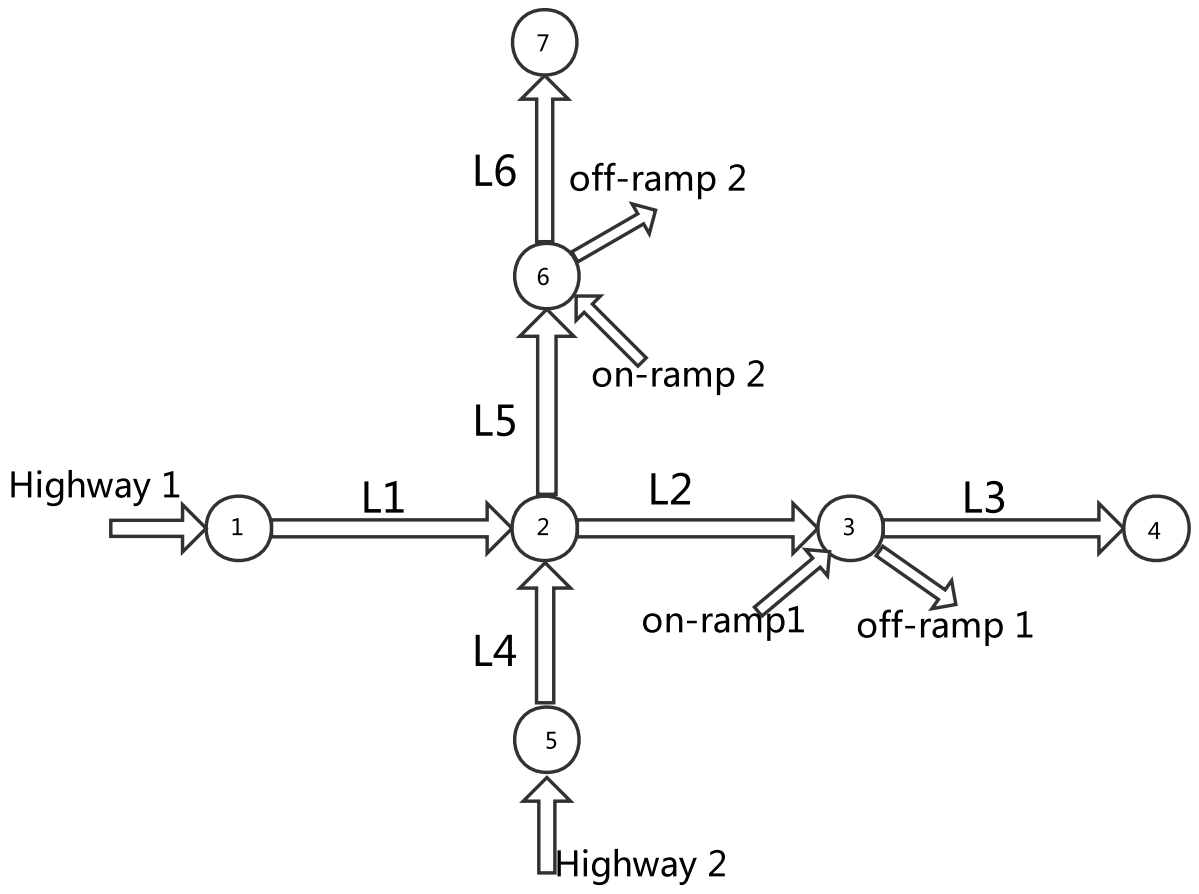

Fig. 2. Freeway network layout.

2. Incoming and outgoing link sets: $L^2_{\text{in}} = \{1, 4\}$, $L^2_{\text{out}} = \{2, 5\}$; $L^3_{\text{in}} = \{2, 7\}$, $L^3_{\text{out}} = \{3, 9\}$; $L^6_{\text{in}} = \{5, 8\}$, $L^6_{\text{out}} = \{6, 10\}$; note that links 7, 8, 9 and 10 are on-ramp 1, on-ramp 2, off-ramp1 and off-ramp2, respectively;

3. The transition matrices at nodes connecting on-ramps and off-ramps are $\boldsymbol{P}^3 = \boldsymbol{P}^6 = \begin{bmatrix} 0.80 & 1.00 \\ 0.20 & 0.00 \end{bmatrix}$;

4. Covariance of turning ratio matrices: $\boldsymbol{\Gamma}^2 = \begin{bmatrix} 0.005 & 0.001 \\ 0.001 & 0.005 \end{bmatrix}$, $\boldsymbol{\Gamma}^3 = \begin{bmatrix} 0.005 & 0 \\ 0 & 0 \end{bmatrix}$, $\boldsymbol{\Gamma}^6 = \begin{bmatrix} 0.005 & 0 \\ 0 & 0 \end{bmatrix}$;

5. The confidence level for each SOC constraint is equal.

Since we assume that the vehicles from an on-ramp will not depart the freeway from the off-ramp at the same node, the covariance matrices are sparse at such nodes. The optimization model is as follows

$$\begin{aligned} \min \quad & -\sum_{i=1}^{n_{max}}((\sum_{j \in L} q_{\text{out}}(i,j) + \sum_{j \in L_{\text{in}}} q_{\text{in}}(i,j) \\ & + m(q_{on}(i,1) + q_{on}(i,2)))(n_{max} - i + 1) - hy(i)) \\ s.t. \quad & y(i) \ge n_{lane}(4) q_{\text{out}}(i,1) - n_{lane}(1) q_{\text{out}}(i,4), \quad \forall i \\ & y(i) \ge n_{lane}(1) q_{\text{out}}(i,4) - n_{lane}(4) q_{\text{out}}(i,1), \quad \forall i \\ & q_{on}(i,1) \ge q_{\text{out}}(i,2)/n_{lane}(2), \quad \forall i \in N \\ & q_{on}(i,2) \ge q_{\text{out}}(i,4)/n_{lane}(4), \quad \forall i \in N \\ & q_{\text{out}}(i,3) \le \psi'(\rho_3), \quad \forall i \in N \\ & q_{\text{out}}(i,6) \le \psi'(\rho_6), \quad \forall i \in N \\ & (15)-(17), \quad \forall j \in L_{\text{in}} \\ & (40), (47), (57) \quad \forall j \in L/L_{\text{in}} \\ & q_{\text{out}}(i,j) \ge 0, \quad q_{\text{in}}(i,j) \ge 0 \quad \forall i, \quad j \end{aligned} \tag{37}$$

The first term of the objective function is to maximize the sum of inflows and outflows of all links and the inflows of on-ramps in this network, and the weights $n_{max} - i + 1$ avoid

$$\boldsymbol{\Gamma}_1^p := T^2 \begin{bmatrix} \boldsymbol{\Gamma}(1,1)_{p\times p} & \boldsymbol{\Gamma}(1,2)_{p\times p} & \cdots & \boldsymbol{\Gamma}(1,n_z)_{p\times p} & \boldsymbol{0}_{p\times 1} \\ \boldsymbol{\Gamma}(2,1)_{p\times p} & \boldsymbol{\Gamma}(2,2)_{p\times p} & \cdots & \boldsymbol{\Gamma}(2,n_z)_{p\times p} & \boldsymbol{0}_{p\times 1} \\ \vdots & \vdots & \vdots & \vdots & \vdots \\ \boldsymbol{\Gamma}(n_z,1)_{p\times p} & \boldsymbol{\Gamma}(n_z,2)_{p\times p} & \cdots & \boldsymbol{\Gamma}(n_z,n_z)_{p\times p} & \boldsymbol{0}_{p\times 1} \\ \boldsymbol{0}_{1\times p} & \boldsymbol{0}_{1\times p} & \cdots & \boldsymbol{0}_{1\times p} & \boldsymbol{0}_{1\times 1} \end{bmatrix} \tag{36}$$

the unnecessary stops, i.e. vehicles will move forward as long as it is not completely congested downstream. $m < 1$ implies the vehicles on the freeway have a priority over the vehicles from on-ramps at nodes 3 and 6. We let $m = 0.1$ in this article. We assume the inflows of incoming boundary links and on-ramps are controllable. Although the outflows are also decision variables, we do not really "control" them since they only need to satisfy the compatible constrains. The second term in the objective function combined with the first two constraints is to add a penalty if the outflows of links 1 and 4 are not proportional to their capacity. A large $h$ is required to add a big penalty term for the violation of this condition, and $h = 100$ is used in this example. The third and fourth constraints set the lower bound of the on-ramp inflows as a function of the outflows of the incoming links at the same node. The fifth and sixth constraints set the supply, i.e. the number of vehicles that a node can accommodate during a unit time, of nodes downstream which are nodes 4 and 7 in this case study, and

$$\psi'(\rho) = \begin{cases} C & \rho \in [0, \rho_c] \\ \psi(\rho) & \rho \in [\rho_c, \rho_m], \end{cases} \quad (38)$$

where $\psi'(\rho)$ is the supply function corresponding to the fundamental diagram Equation (2) and $C$ is the capacity. The seventh and eighth constraints indicate that the incoming boundary links satisfy the deterministic compatibility conditions, and other links satisfy the distributionally robust chance constraints.

### B. Results

We studied three different scenarios in terms of the level of service of this network. The first scenario is that the network is under free flow condition; the second scenario is that the network is congested; the third case is that the network is partially congested.

*1) Free Flow Network:* Let the initial densities of every link be equal to their critical densities, and $\boldsymbol{\rho} = 0.8\boldsymbol{\rho_c}$, $\boldsymbol{P}^2 = \begin{bmatrix} 0.80 & 0.27 \\ 0.20 & 0.73 \end{bmatrix}$. This transition matrix makes the inflows of links 2 and 5 proportional to their capacity. Then, we investigate the influence of the uncertainty of turning ratio at one node on other control results. Figure 3 shows the control inputs considering the uncertainties of node 6. The base case is the result without considering the uncertainties, so it is a non-robust control result.

Every 25 points on the horizontal axis indicate the control inputs for an incoming link for 25 time steps. For example, the first 25 points represent the optimal inflow of link 1 during the simulation, and the range of 50-75 indicates the optimal inflow of on-ramp 1. The base case is the results without uncertainties.

The change of the inflows from on-ramp 2 shows that the distributionally robust chance constraints require the original constraints to hold with a high probability for the worst distribution, and this will lead to a more conservative optimal solution. Due to the maximization objective function, the sum of the incoming flows from link 5 and the inflow of on-ramp 2 should be equal to the supply of link 6. If we are not certain about the turning ratio at node 6, we have to decrease the on-ramp inflows to ensure link 6 can accommodate all the incoming vehicles even for the worst case. This is the reason why the inflow of on-ramp 2 is lower than the deterministic model at the beginning of the simulation. For the same reason, the constraints (47) restrict the outflow of link 6 to make it less than the Moskowitz solution from the upstream conditions. Therefore, its outflow is less than its capacity although it is in free flow condition downstream. This constraint adds a shockwave moving backward, and the supply of link 6 will decrease once the shockwave reaches its upstream end. This is the reason why there is a drop on the inflows of on-ramp 2 appearing in the middle of the simulation. In addition, under free flow condition, it is shown that there is little impact of an intersection (node 6) on upstream links (links 1 and 4) or other freeways (on-ramp 1). Note that the drop of inflows of links 1 and 4, appearing at the end of the simulation, results from the deviation of the transition equation. In this case, $P^2_{11}C_1 + P^2_{12}C_4$ is a little larger than $C_2$; this will block a small part of vehicles and induces a shockwave on both links 1 and 4.

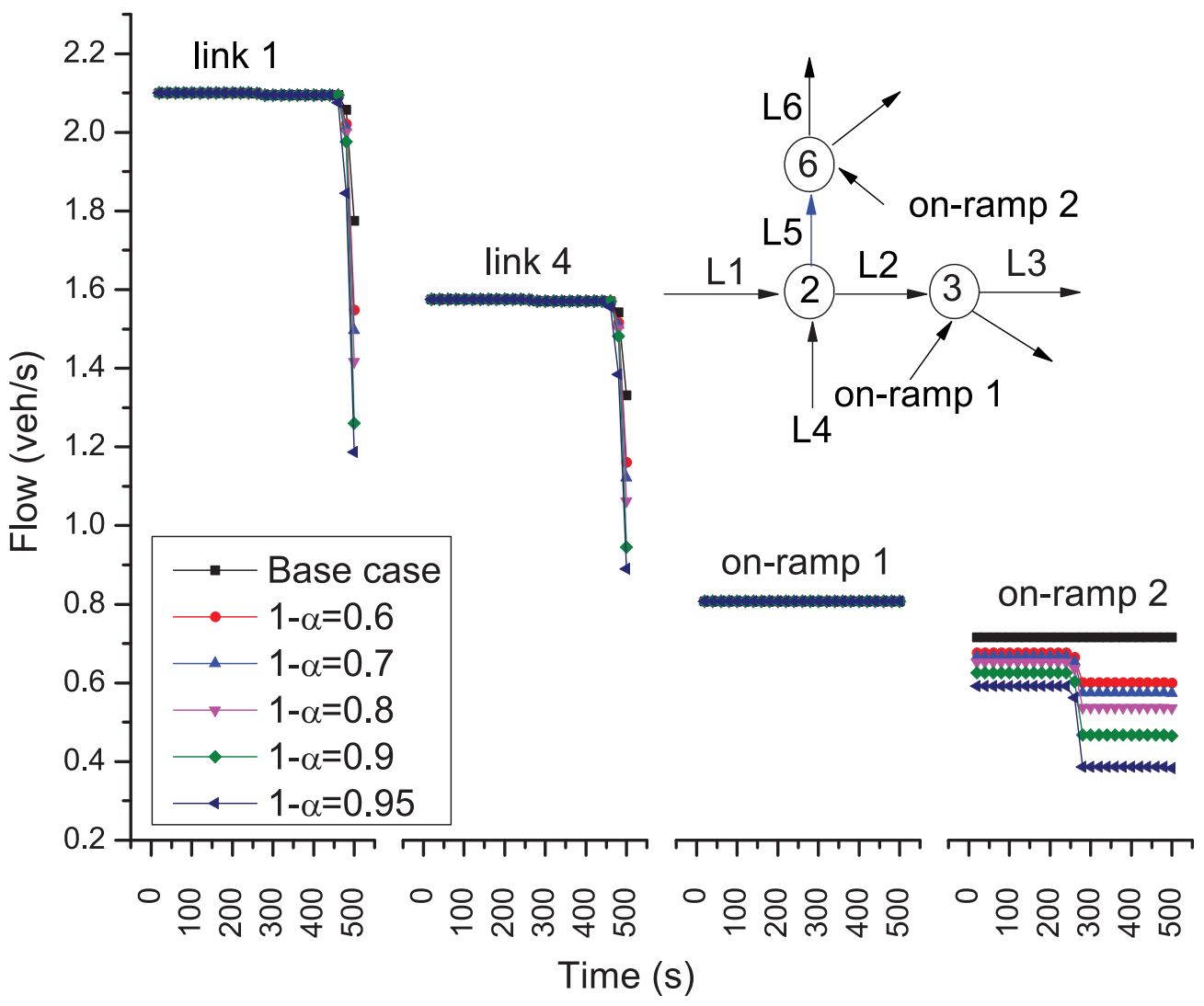

Fig. 3. Optimal inflows considering turning ratio uncertainty at node 6.

The control inputs considering the uncertainty at node 2 is shown in Figure 4. Due to the same reason mentioned above, the optimal inflows of links 1 and 4 drop at some point. Consequently, the inflows of two on-ramps increases since fewer vehicles upstream merge to the node. Therefore, the uncertainties decrease the inflow of the incoming links and increase the inflows of on-ramps downstream.

*2) Congested Network:* Let us assume a congested network with $\boldsymbol{\rho} = 4\boldsymbol{\rho_c}$, $\boldsymbol{P}^2 = \begin{bmatrix} 0.80 & 0.27 \\ 0.20 & 0.73 \end{bmatrix}$. The corresponding control inputs are shown in Figures 5 and 6. Figure 5 shows a similar phenomenon as Figure 3. Unlike the free flow case, Figure 6 shows that the on-ramps downstream are not impacted by the intersection upstream for the congested case. This is because although the inflows of links 1 and 4 are reduced, the whole network is still congested which means there are enough vehicles from the freeway links (links 2 and 5) merging with on-ramps during the simulation period. If the simulation horizon increases to some extent, the inflows of links 1 and 4 will further decrease and the inflows of on-ramps will be expected to increase.

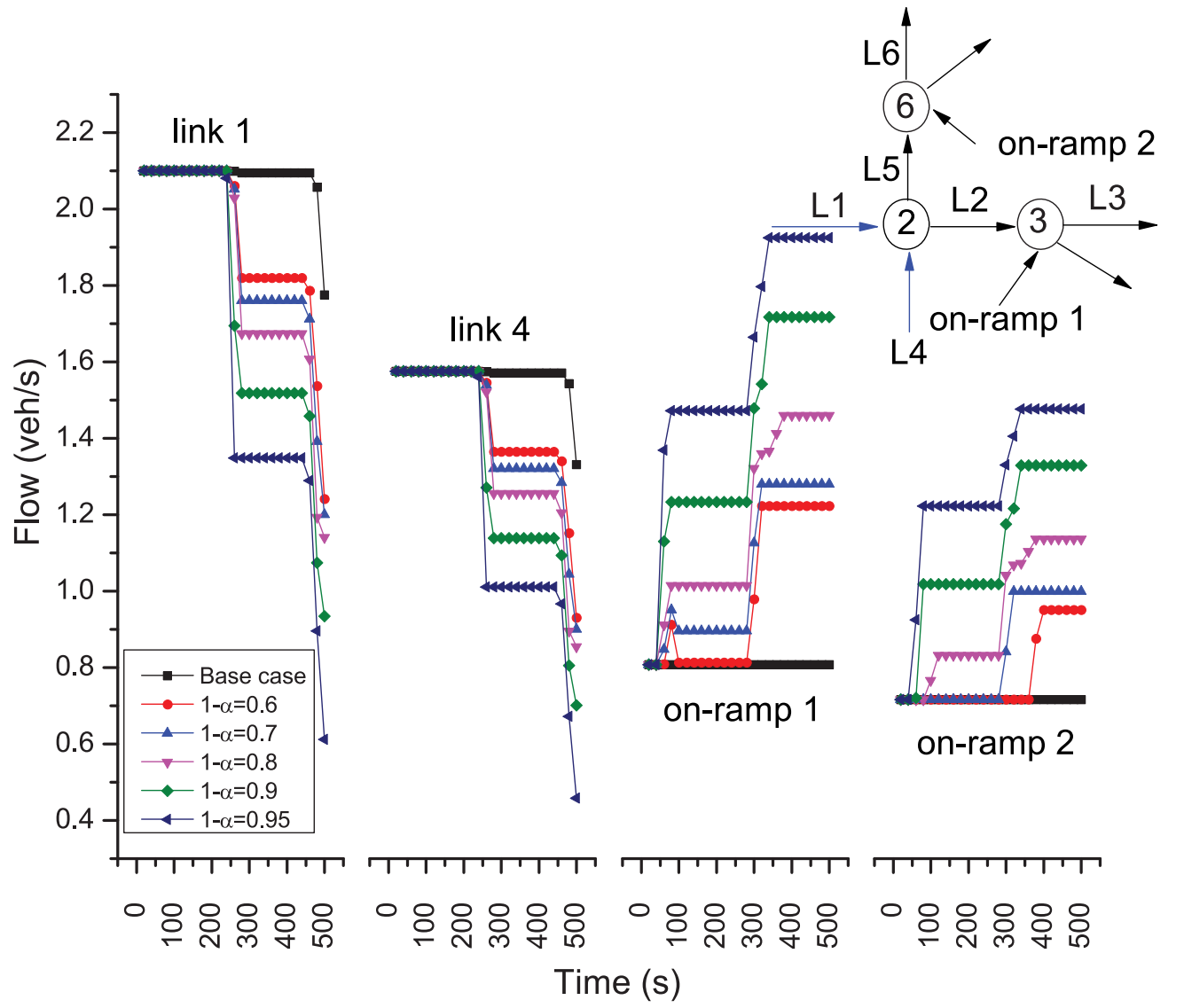

Fig. 4. Optimal inflows considering turning ratio uncertainty at node 2.

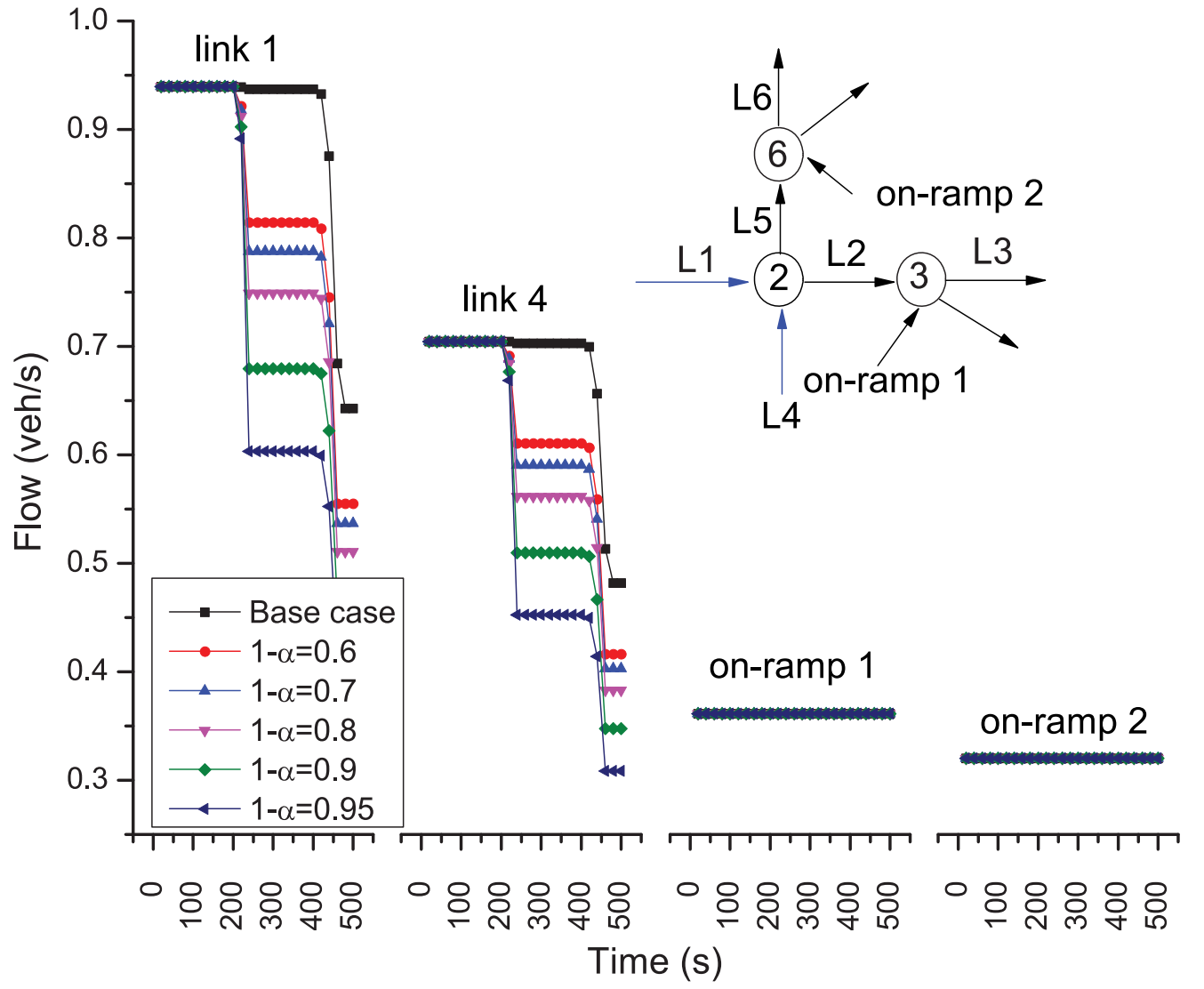

Fig. 6. Optimal inflows considering turning ratio uncertainty at node 2.

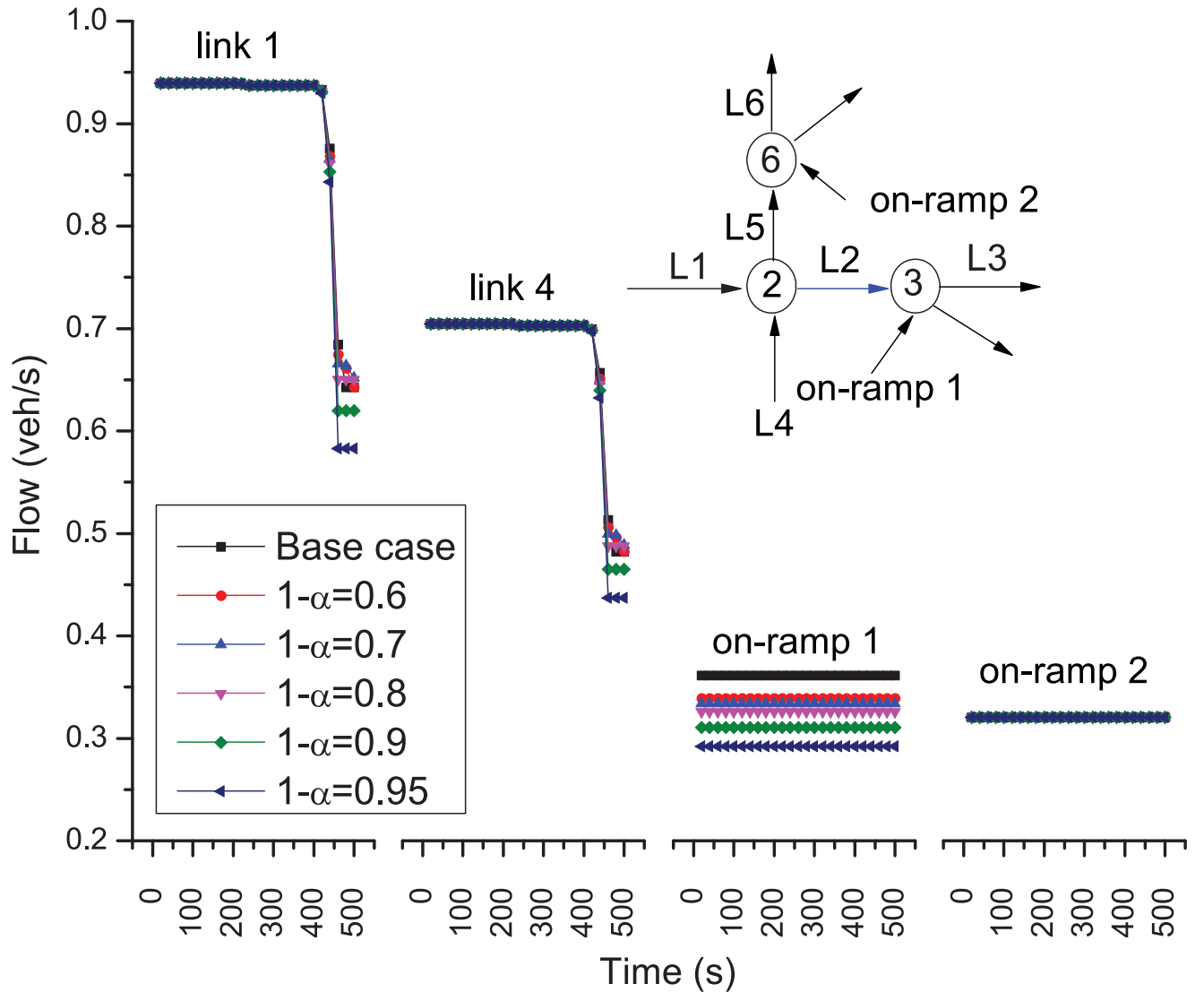

Fig. 5. Optimal inflows considering turning ratio uncertainty at node 3.

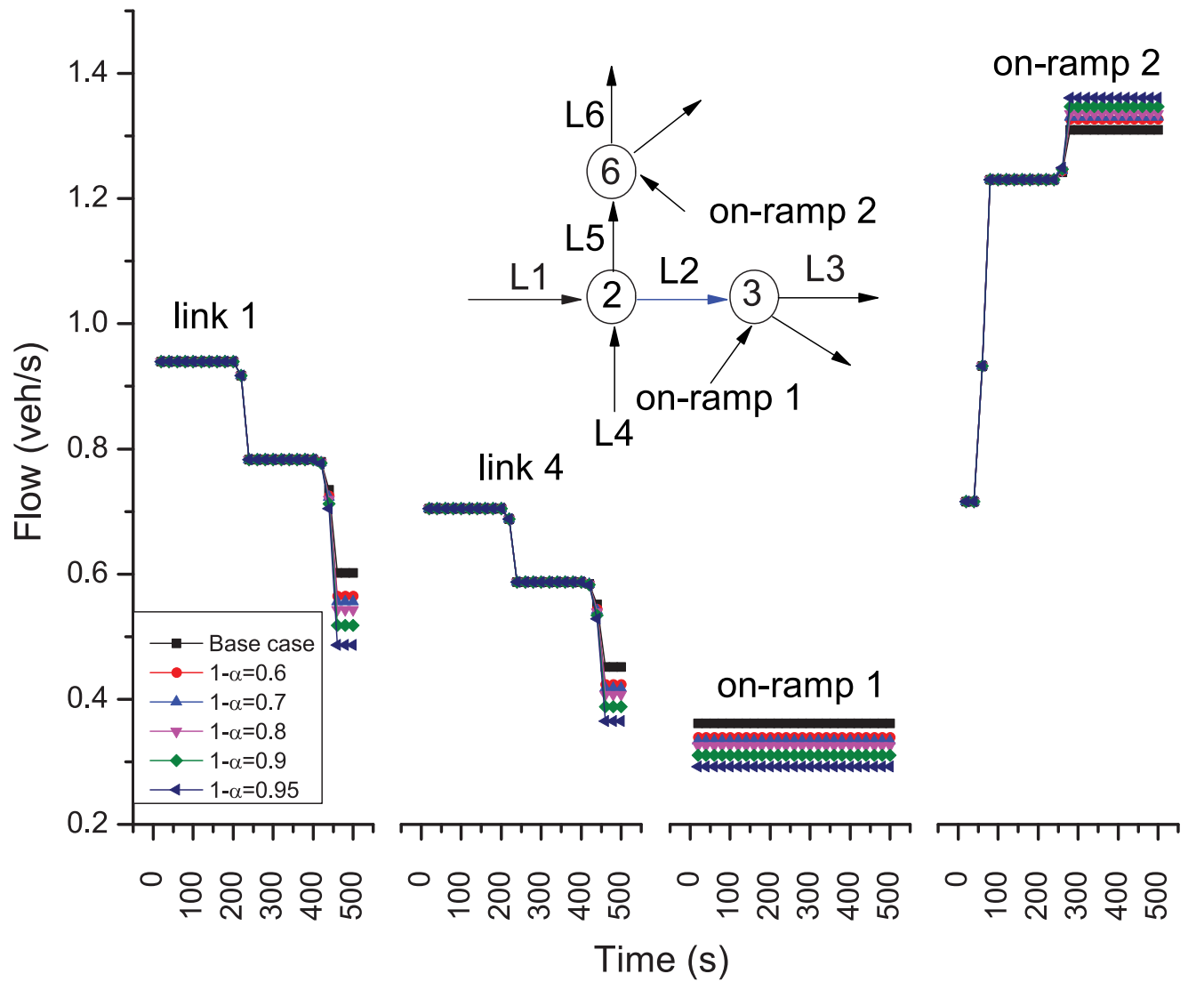

Fig. 7. Optimal inflows considering turning ratio uncertainty at node 3.

*3) Partially Congested Network:* Let us consider another network with

$$\rho_l = \begin{cases} 4\rho_c & l \in \{1,2,3,4\} \\ 0.8\rho_c & l \in \{5,6\}, \end{cases} \qquad \boldsymbol{P}^2 = \begin{bmatrix} 0.60 & 0.80 \\ 0.40 & 0.20. \end{bmatrix}$$

In this scenario, links 5 and 6 are under free flow condition, and all other links are congested since the demand of link 2 is higher than its supply, which is implied by the transition matrix at node 2. Figure 7 shows the control inputs with random turning ratios at node 3. This uncertainty reduces the on-ramp flows at the same node and increases the flows of on-ramp 2. Unlike the completely congested network, link 2 in this case blocks the vehicles on links 1 and 4, and since links 5 and 6 are under free flow conditions, on-ramp 2 can send more vehicles to maximize the total throughput. The two drops on links 1 and 4 are subjected different regimes or shockwaves. At the beginning of the simulation, shockwaves are generated at both nodes 2 and 3. The shockwave at node 2 originates from the fact that, due to the transition matrix and the congestion condition, the number of vehicles passing through node 2 on links 1 and 4 is less than the number of vehicles desiring to take that route. The shockwave at node 3 is generated by the same reason in the base case. For the robust model, robust constraints are another reason of this shockwave. Both shockwaves from nodes 2 and 3 move backward and will arrive the upstream end of boundary links (links 1 and 4) and restrict their inflows. Therefore, the inflows decrease by the same scale the first time while the second reduction increases with the confidence level $1-\alpha$.

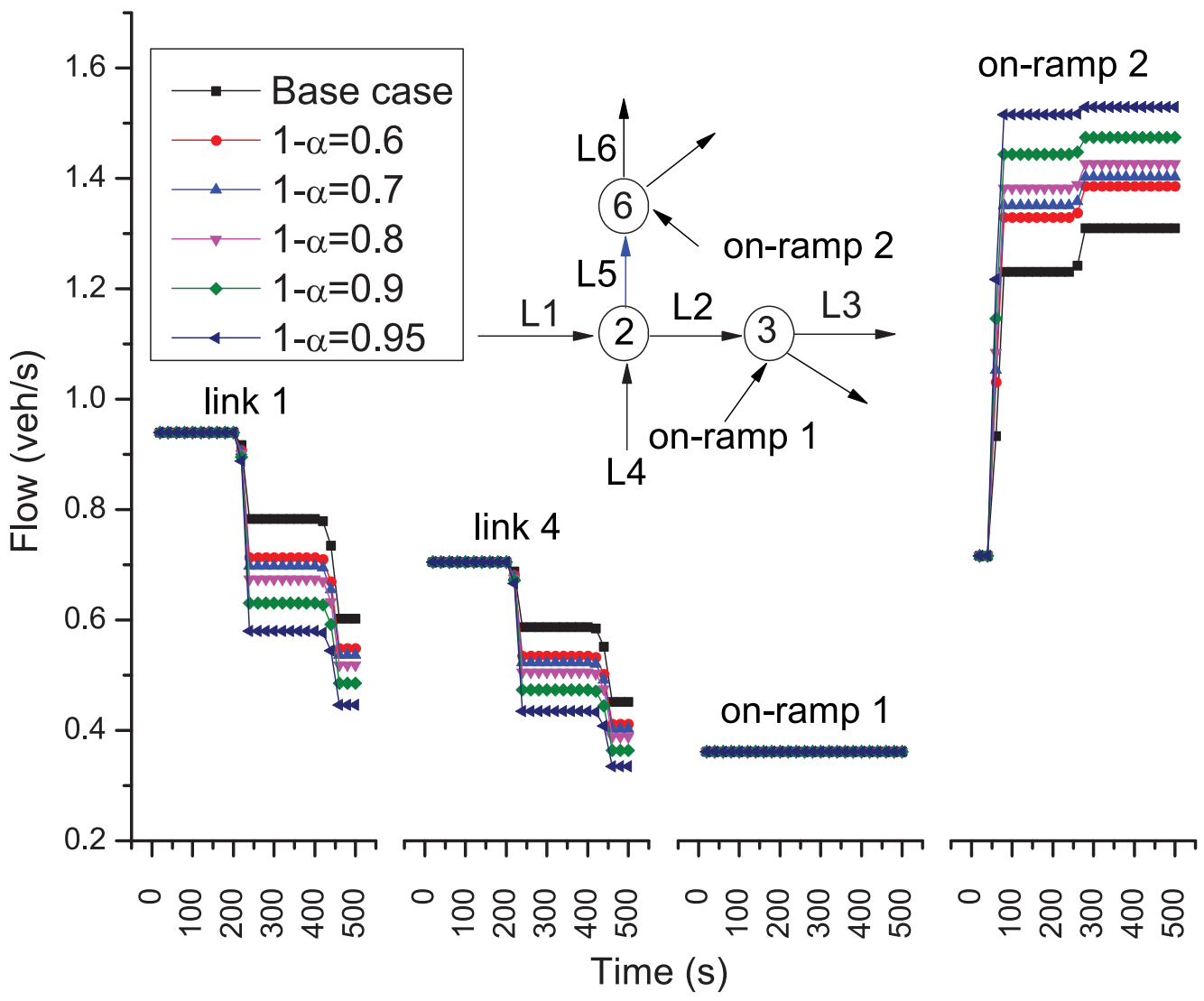

Fig. 8. Optimal inflows considering turning ratio uncertainty at node 2.

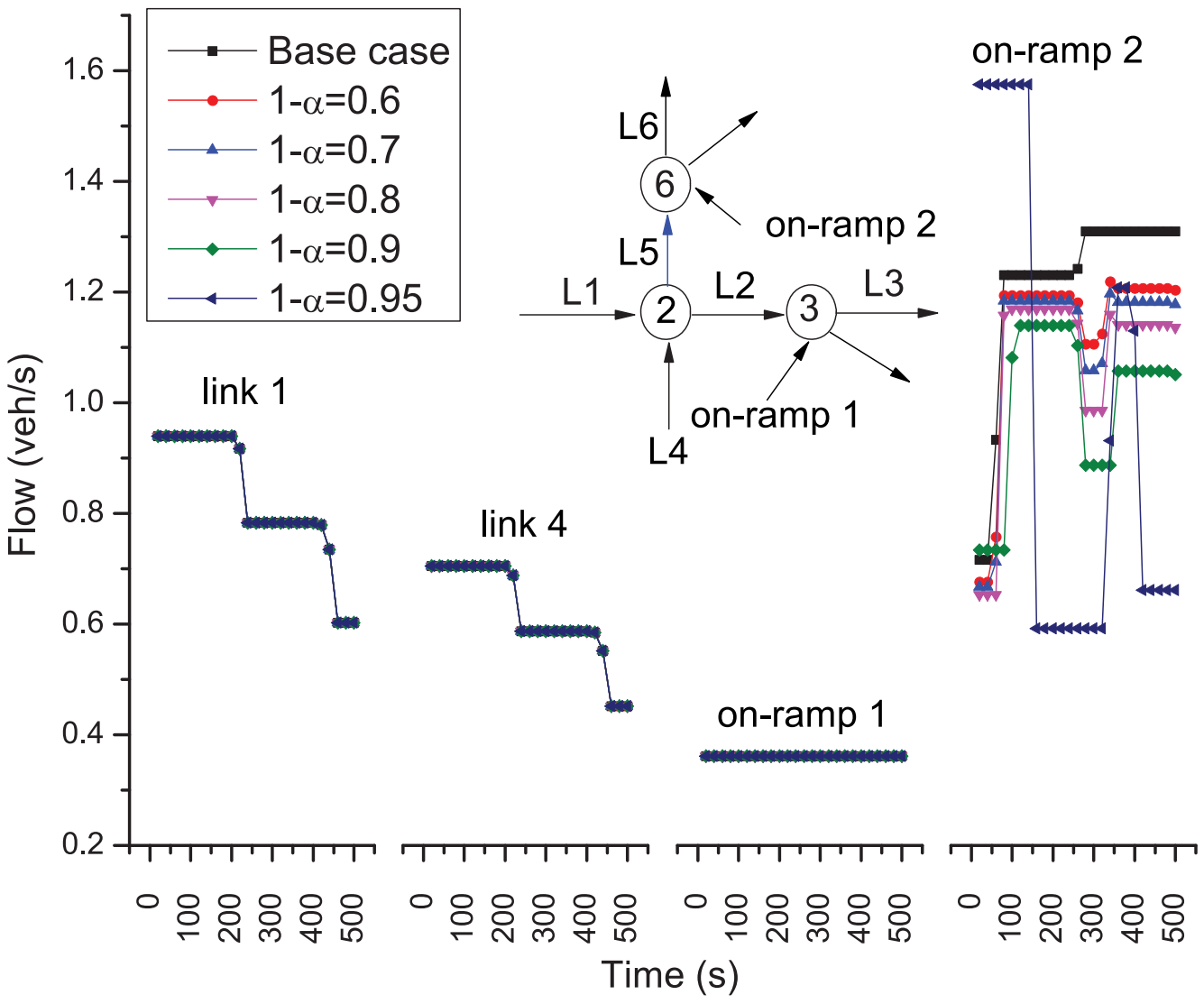

Fig. 9. Optimal inflows considering turning ratio uncertainty at node 6.

Figure 8 shows the control inputs with random turning ratios at node 2. Similarly, the on-ramp inflows downstream increases if it is under free flow condition and would not be impacted if it is congested.

Figure 9 shows the control inputs considering the uncertainties at node 6. Two phenomena need to be explained:

1. When the confidence level is high, e.g. $1-\alpha = 0.95$, the inflows of on-ramp 2 at the beginning is equal to the supply of link 6, and the outflow of link 5 is zero. This is because that link 5 is under free flow condition, and the number of vehicles on it at the beginning of simulation is low. Therefore, for a fixed and low outflow of link 5, due to the constraint (47), the outflow of link 6 is very small if $1-\alpha$ is large. If we block the link 5 and proceed the vehicles from on-ramp 2, this will increase the outflows of link 6 in first several time steps substantially since we assume all the vehicle from on-ramp 2 will flow onto link 6. When this contribution to the objective value is more significant than the loss from blocking link 5, it is reasonable to execute the corresponding decision.

2. When $1-\alpha$ is small, there is a reduction on the inflow curve of on-ramp 2 followed by an increase. This is also caused by two distinct shockwaves. The first shockwave is generated when the vehicles on link 5 at the beginning of the simulation arrive at the downstream node of link 6. Let us call this shockwave $S1$ and the equivalent density downstream $\mu_1$. In the following, since links 1 and 4 are blocked, there are not enough vehicles to come onto link 6 from link 5, so more vehicles from on-ramp 2 can flow in. Again, since there are no uncertainties in the turning ratios from on-ramp 2, the outflow of link 6 increases if more vehicles merge onto link 6. Therefore, a short time after the first shockwave is generated, a second shockwave $S2$ is formed, and the corresponding density downstream $\mu_2 < \mu_1$. When these shockwaves touch the upstream node 6, the inflow of on-ramp 2 changes correspondingly. Since $\mu_2 < \mu_1$, the inflows should decrease and then increase.

### C. Performance Demonstration

To demonstrate the performance of the proposed model, the robust control and non-robust control results are implemented on a scenario in which the real turning ratio deviates from the mean matrix, and relevant traffic states are compared. To have a fair comparison, instead of using the Moskowitz solution, we employ the cell transmission model (CTM) [45], [46], which has been widely recognized as a benchmark simulation [8], [18], to simulate traffic evolution. In the CTM model, the time step is 5 s and the cell length is 200 m. All other model parameters are defined in Section IV-A.

Let us consider the scenario shown in Figure 8 in which the turning ratios at node 2 are random. The non-robust control, i.e., the base case, and robust control with confidence level of 0.9 shown in Figure 8 are implemented on a case where the realization of the turning ratio matrix is $\boldsymbol{P}^2 = \begin{bmatrix} 0.70 & 0.85 \\ 0.30 & 0.15 \end{bmatrix}$ while the mean matrix is $\boldsymbol{P}^2 = \begin{bmatrix} 0.60 & 0.80 \\ 0.40 & 0.20 \end{bmatrix}$. In this case, this matrix sends more vehicles to link 2 than its capacity, so the non-robust control may cause congestion that continues moving backward. Figure 10 shows the number of vehicles that blocked at the upstream nodes of links 1 and 4, which denotes the excessive number of vehicles from the optimal solutions that cannot be accommodated by links 1 and 4. Compared to the non-robust solution, the robust control generates lower inflows for links 1 and 4 by considering the risk of congestion brought by the uncertainties at node 2, and this feature avoids the risk of blocking vehicles at the entry nodes.

Apart from the benefit of reducing the probability of congestion, the robust control can also increase the optimal value, i.e., total throughput, at the same time. Figure 11 shows the comparison of throughput of link 6. Although the non-robust control sends more vehicle through links 1 and 4, not enough vehicles can arrive onto link 5 because: 1. the real turning ratio matrix sends more vehicles to link 2; 2. node 2 is congested and prevents some vehicles coming through. On the contrary, the robust control solution is capable of detecting such uncertainty and increases the on-ramp flow, shown in Figure 8, to maximize the objective function.

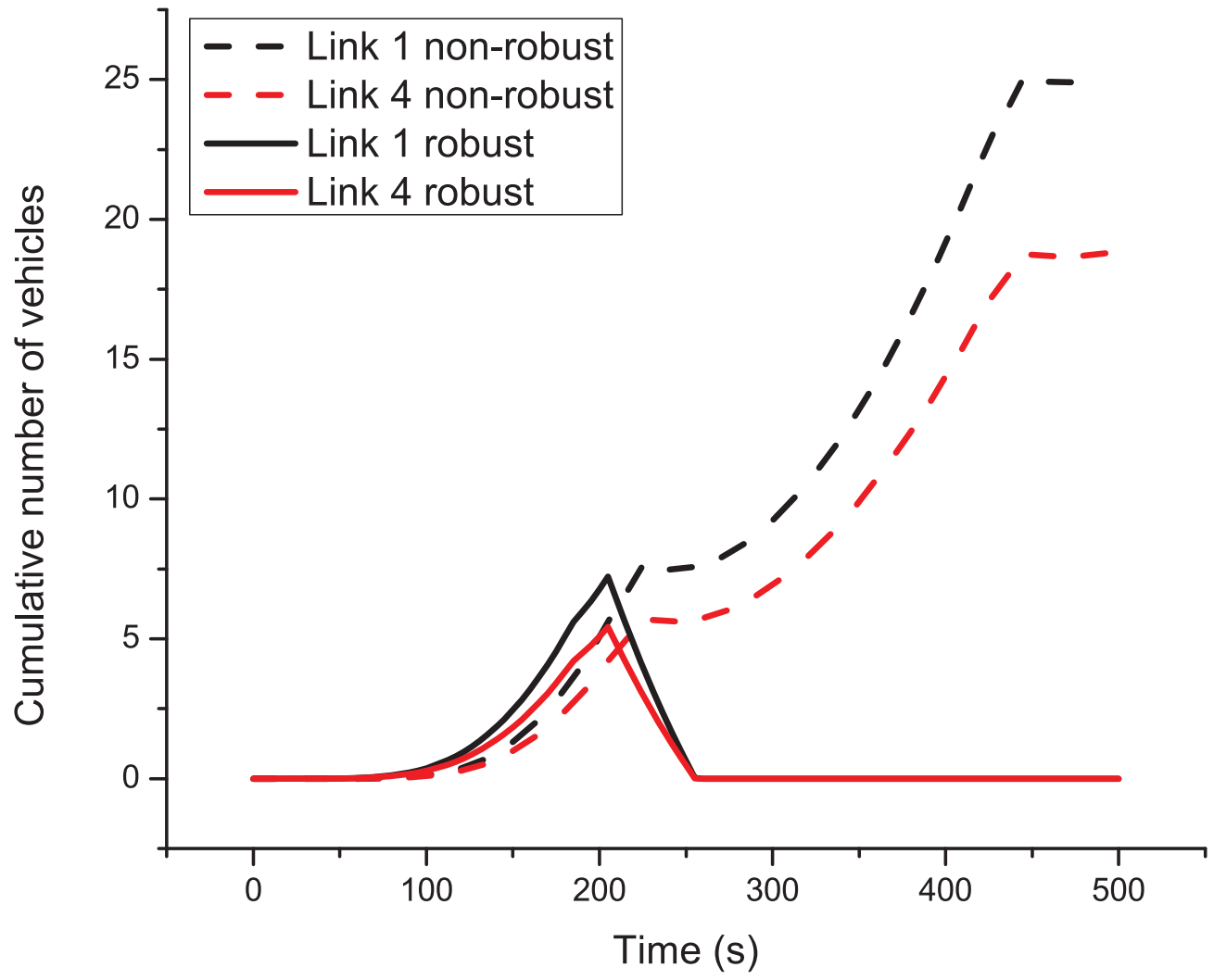

Fig. 10. Number of vehicles blocked at entry nodes.

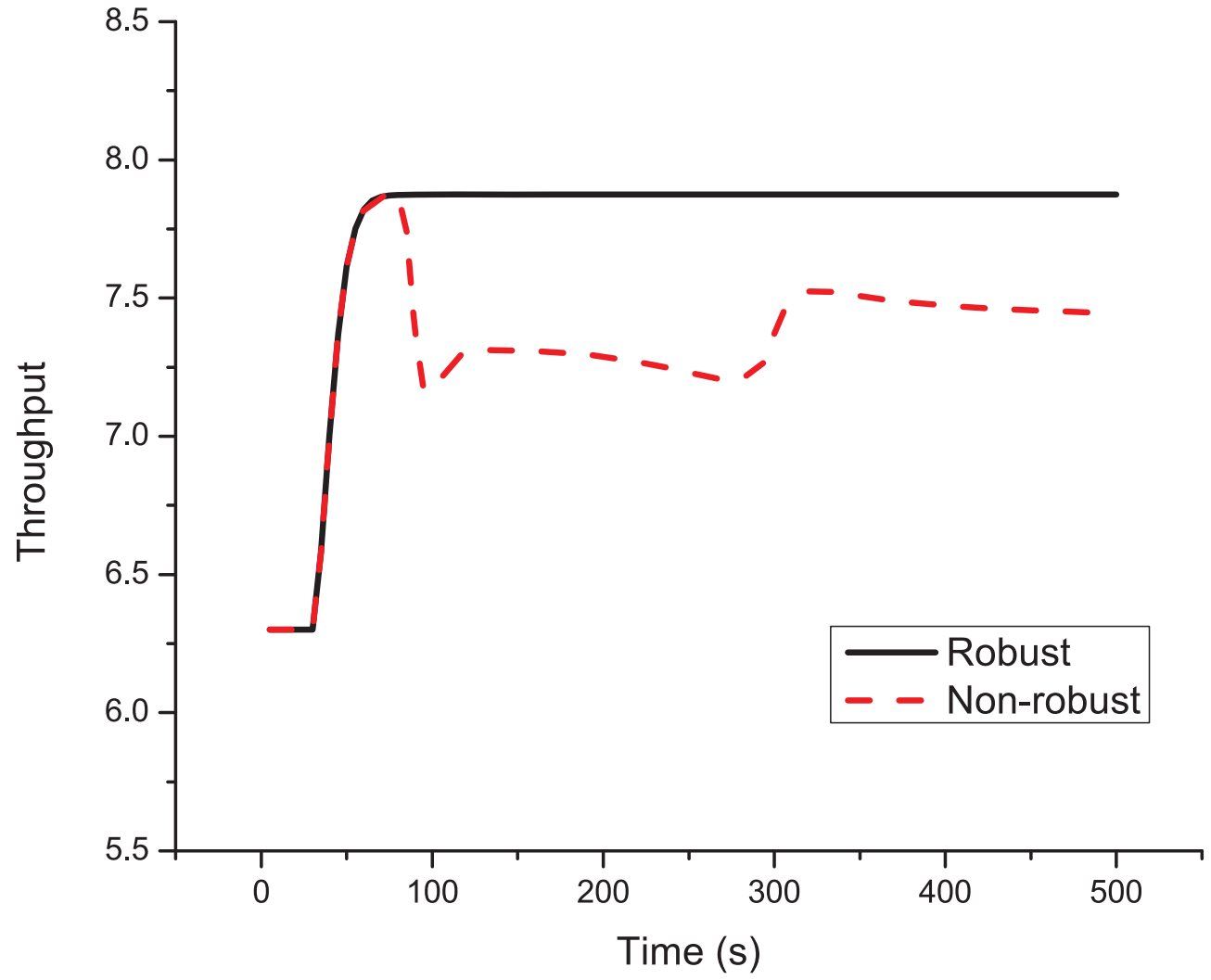

Fig. 11. Throughput of link 6.

Overall, the robust control model outperforms the non-robust control in terms of reducing congestion and improving throughput.

## V. Case Study on an Urban Network

In addition to freeway networks, this section applies the proposed framework on an urban network and the corresponding inflow controls are studied.

### A. Network and Problem Settings

A sub-network of downtown Austin, TX, shown as the blue square in Figure 12, is employed. This network consists of 55 links and 20 nodes, and all the nodes are signalized. The link attributes, such as length and number of lanes, and signal timings, such as splits and offsets, are obtained from the database of VISTA [47] administered by the University of Texas at Austin. In this network, except for Congress Avenue, all streets are one-way streets. The model parameters are as follows: free flow speed $v_f = 13.5$ m/s, capacity $C = 0.3375$ veh/s/lane, critical density $\rho_c = 0.025$ veh/m/lane, jam density $\rho_m = 0.125$ veh/m/lane and congestion speed $w = -3.86$ m/s. For simplicity, let the link lengths be unanimous and equal to 128 m. We divide each link into 2 even segments. Let the initial densities be 0 and $\rho_c$ in the upstream and downstream segments, respectively. Additionally, the simulation time is 300 s and is divided into 75 even time steps.

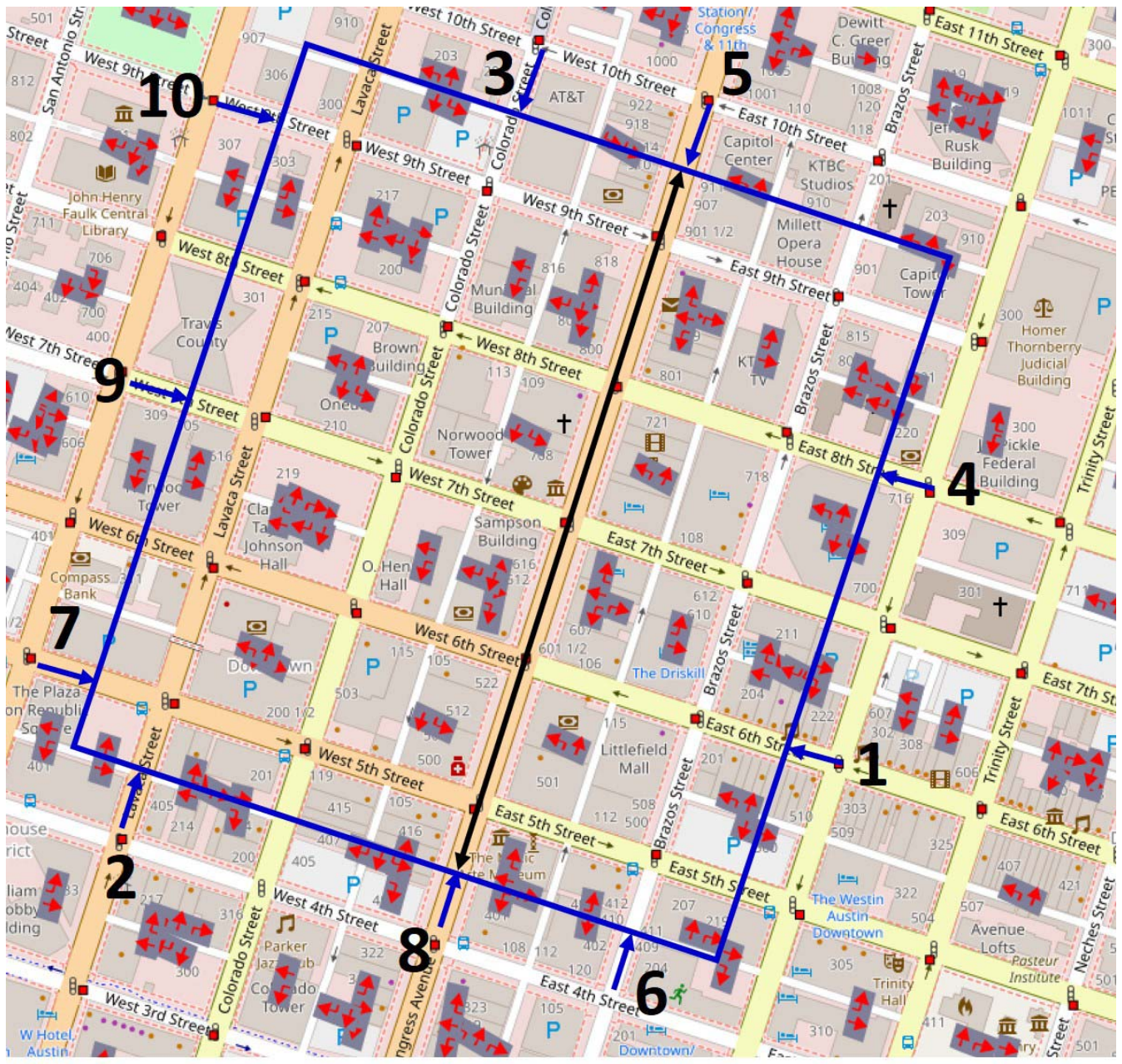

Fig. 12. A sub-network of downtown Austin.

Similarly as the freeway network used in Section IV, we assume the inflows of incoming boundary links, shown as the blue arrows in Figure 12, are continuous and controllable. The optimization model is

$$\begin{aligned}
\min \quad & -\sum_{i=1}^{n_{max}}((\sum_{j\in L} q_{\text{out}}(i,j) + \sum_{j\in L_{\text{in}}} q_{\text{in}}(i,j))(n_{max} - i + 1)) \\
& + \omega \sum_{i=1}^{n_{max}} \sum_{j\in L_{\text{in}}} (q_{\text{in}}(i,j) - q_{\text{in}}(i+1,j))^2 \\
s.t. \quad & q_{\text{out}}(i,j) \le e_j C_j, \quad \forall i \in N, j \in L_{\text{out}} \\
& (15) - (17), \quad \forall j \in L_{\text{in}} \\
& (40), (47), (57) \quad \forall j \in L/L_{\text{in}} \\
& q_{\text{out}}(i,j) \le C_j s(i,j), \quad \forall i \in N, j \in L \\
& q_{\text{out}}(i,j) \ge 0, \quad q_{\text{in}}(i,j) \ge 0 \quad \forall i, \quad j
\end{aligned} \tag{39}$$

The first term of the objective function is to maximize the sum of weighted outflows and inflows; the second term is a quadratic function which is used to reduce the inflow fluctuation, i.e. smooth the inflows. $\omega$ is the corresponding weight.

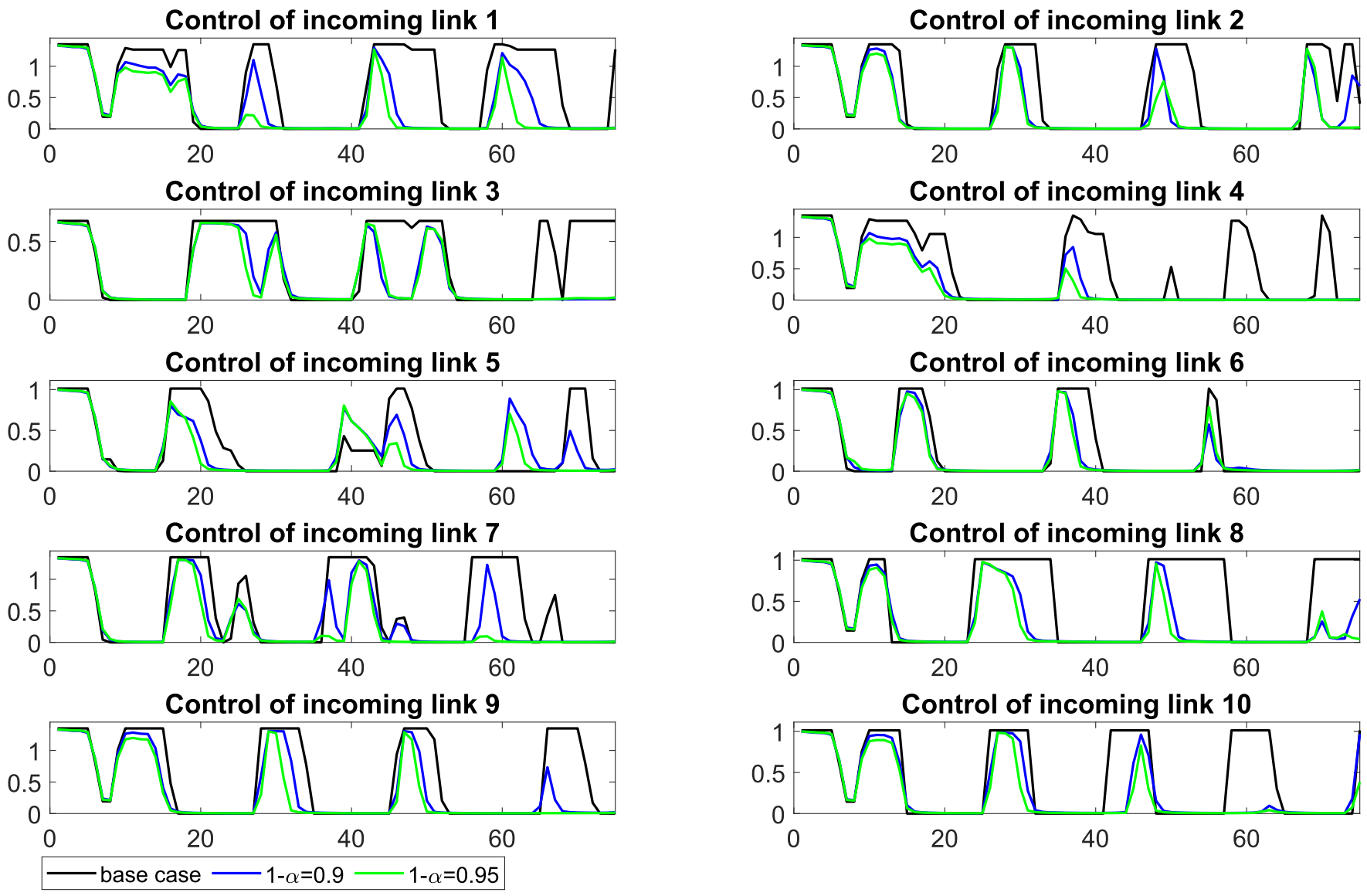

Fig. 13. Control inputs with $\omega = 0.2$, showing flow (veh/s) at each time index.

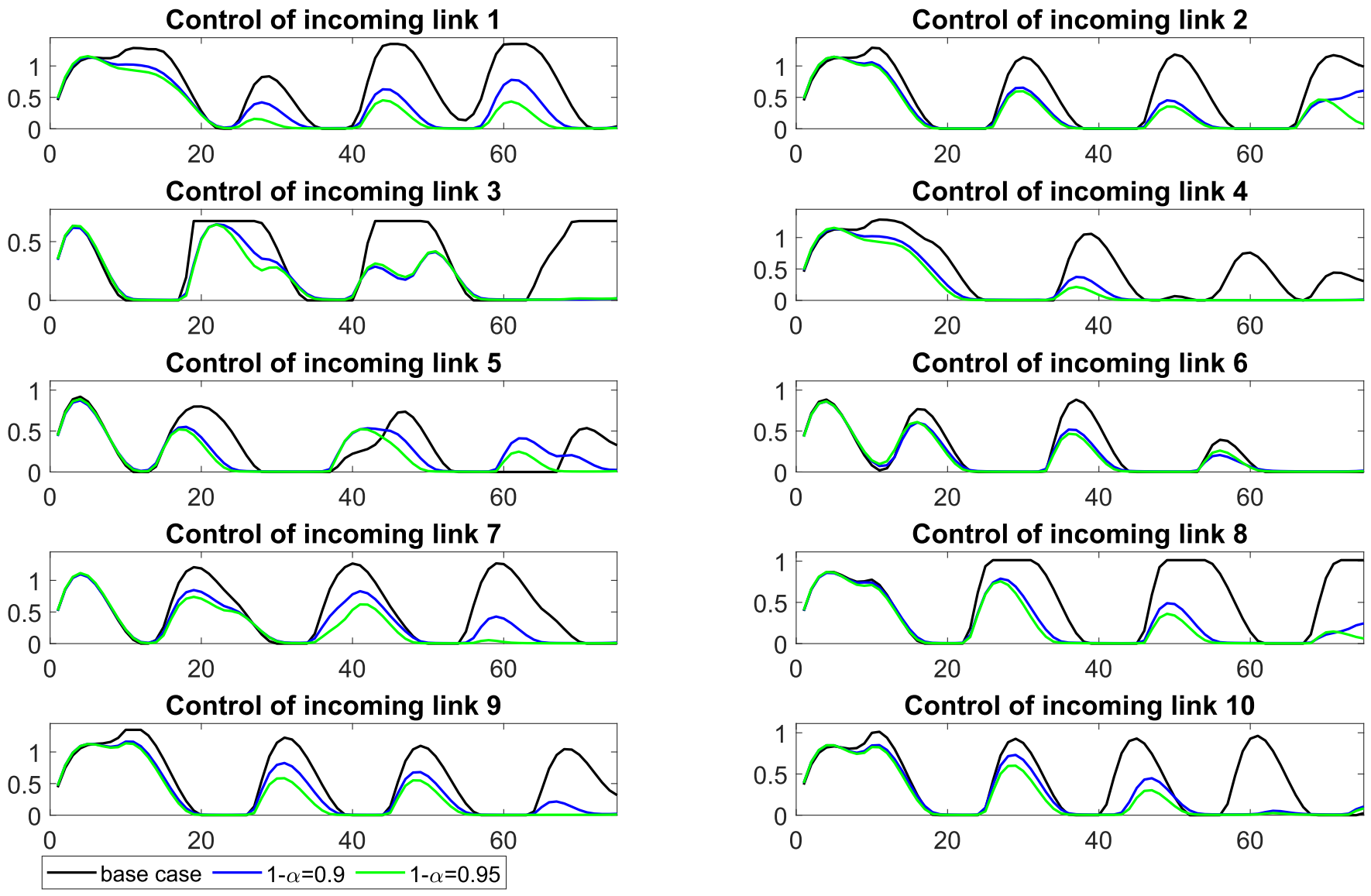

Fig. 14. Control inputs with $\omega = 20$, showing flow (veh/s) at each time index.

The first constraint adds an upper bound for the outflows of outgoing boundary links as a proportion of their capacities. In this example, $e_j = 0.8$ for all the outgoing boundary links. The fourth constraint adds the signal constraints where $s(i, j)$'s are binary parameters. $s(i, j) = 1$ when the phase serving $j$ at time $i$ is green, and $s(i, j) = 0$ otherwise. The remaining constraints are the same as the freeway network in Section IV. Note that since most of the streets are one-way, all the vehicles from one link are served by the same phase. For a network consisting of two-way streets, we need to group movements that can be served by the same phase and treat the groups as separate links.

We assume all the vehicles go straight at an intersection with a probability of 0.8, and the remaining turning ratios are equal to $\frac{0.2}{|I|-1}$ where $|I|$ is the number of outgoing links at intersection $I$. Note that these assumptions on turning ratios are not necessary for the proposed model, and any other estimation could be used based on users' preferences. In addition, we assume the variance of all $P(i, j)$'s is equal to 0.005, and the covariance between any two incoming links $P(i, j_1)$ and $P(i, j_2)$ is equal to 0.001.

## B. Results

Figures 13 and 14 show the optimal controls on those 10 incoming boundary links from Figure 12, over the simulation time with different weights $\omega$. The horizontal axis is the index of time steps while the vertical axis is the optimal inflows (veh/s). The differences between the ranges over the vertical axis result from the different number of lanes. There

$$\begin{cases} M_{M_k}(0, x_p) \geq M_p(0, x_p) & \forall (k,p) \in K^2 \\ M_{M_k}(pT, \chi) \geq \beta_p(pT, \chi) & \forall k \in K, \quad \forall p \in N \\ M_{M_k}(\frac{\chi - x_k}{v_f}, \chi) \geq \beta_p(\frac{\chi - x_k}{v_f}, \chi) & \forall k \in K, \quad \forall p \in N \\ \qquad\qquad \text{s.t. } \frac{\chi - x_k}{v_f} & \in [(p-1)T, pT] \\ \kappa_\alpha \|(\boldsymbol{\Gamma}_1^p)^{\frac{1}{2}} \tilde{\boldsymbol{x}}_1^p\|_2 + \hat{\boldsymbol{d}}_1^p(\tilde{\boldsymbol{x}}_1^p) \leq 0, & \forall k \in K, \quad \forall p \in N \\ \kappa_\alpha \|(\boldsymbol{\Gamma}_2^p)^{\frac{1}{2}} \tilde{\boldsymbol{x}}_2^p\|_2 + \hat{\boldsymbol{d}}_2^p(\tilde{\boldsymbol{x}}_2^p) \leq 0 & \forall k \in K, \quad \forall p \in N \\ \qquad\qquad \text{s.t. } \frac{\xi - x_{k-1}}{w} & \in [(p-1)T, pT] \end{cases} \tag{40}$$

$$\begin{aligned} \tilde{\boldsymbol{x}}_1^p := & [\overbrace{q_{\text{out}}(1, L_{\text{in}}^z(1)), q_{\text{out}}(2, L_{\text{in}}^z(1)), \ldots, q_{\text{out}}(p, L_{\text{in}}^z(1)),}^{p} \overbrace{\cdots\cdots}^{(n_z-2)\times p}, \\ & \overbrace{q_{\text{out}}(1, L_{\text{in}}^z(n_z)), q_{\text{out}}(2, L_{\text{in}}^z(n_z)), \ldots, q_{\text{out}}(p, L_{\text{in}}^z(n_z))}^{p}, 1]^T, \end{aligned} \tag{41}$$

$$\begin{aligned} \hat{\boldsymbol{d}}_1^p := & [\overbrace{T\boldsymbol{P}^z(L_{\text{out}}^z(i), L_{\text{in}}^z(1)), \ldots, T\boldsymbol{P}^z(L_{\text{out}}^z(i), L_{\text{in}}^z(1)),}^{p} \overbrace{\cdots\cdots}^{(n_z-2)\times p}, \\ & \overbrace{T\boldsymbol{P}^z(L_{\text{out}}^z(i), L_{\text{in}}^z(n_z)), \ldots, T\boldsymbol{P}^z(L_{\text{out}}^z(i), L_{\text{in}}^z(n_z))}^{p}, -M_{M_k}(pT, \xi)], \end{aligned} \tag{42}$$

$$\boldsymbol{\Gamma}_1^p := T^2 \begin{bmatrix} \boldsymbol{\Gamma}(1,1)_{p\times p} & \boldsymbol{\Gamma}(1,2)_{p\times p} & \cdots & \boldsymbol{\Gamma}(1,n_z)_{p\times p} & \boldsymbol{0}_{p\times 1} \\ \boldsymbol{\Gamma}(2,1)_{p\times p} & \boldsymbol{\Gamma}(2,2)_{p\times p} & \cdots & \boldsymbol{\Gamma}(2,n_z)_{p\times p} & \boldsymbol{0}_{p\times 1} \\ \vdots & \vdots & \vdots & \vdots & \vdots \\ \boldsymbol{\Gamma}(n_z,1)_{p\times p} & \boldsymbol{\Gamma}(n_z,2)_{p\times p} & \cdots & \boldsymbol{\Gamma}(n_z,n_z)_{p\times p} & \boldsymbol{0}_{p\times 1} \\ \boldsymbol{0}_{1\times p} & \boldsymbol{0}_{1\times p} & \cdots & \boldsymbol{0}_{1\times p} & \boldsymbol{0}_{1\times 1} \end{bmatrix}, \tag{43}$$

$$\tilde{\boldsymbol{x}}_2^p := \tilde{\boldsymbol{x}}_1^p, \tag{44}$$

$$\begin{aligned} \hat{\boldsymbol{d}}_1^p := & [\overbrace{T\boldsymbol{P}^z(L_{\text{out}}^z(i), L_{\text{in}}^z(1)), \ldots,}^{p-1} \boldsymbol{P}^z(L_{\text{out}}^z(i), L_{\text{in}}^z(1))t_1, \overbrace{\cdots\cdots}^{(n_z-2)\times p}, \\ & \overbrace{T\boldsymbol{P}^z(L_{\text{out}}^z(i), L_{\text{in}}^z(n_z)), \ldots,}^{p-1} \boldsymbol{P}^z(L_{\text{out}}^z(i), L_{\text{in}}^z(n_z))t_1, -M_{M_k}(pT, \xi)], \end{aligned} \tag{45}$$

$$\boldsymbol{\Gamma}_2^p := \begin{bmatrix} T^2\boldsymbol{\Gamma}(1,1)_{(p-1)\times(p-1)} & (Tt_1\boldsymbol{\Gamma}(1,1))_{(p-1)\times 1} & T^2\boldsymbol{\Gamma}(1,2)_{(p-1)\times(p-1)} & (Tt_1\boldsymbol{\Gamma}(1,2))_{(p-1)\times 1} & \cdots & \boldsymbol{0}_{(p-1)\times 1} \\ (Tt_1\boldsymbol{\Gamma}(1,1))_{1\times(p-1)} & (t_1^2\boldsymbol{\Gamma}(1,1))_{1\times 1} & (Tt_1\boldsymbol{\Gamma}(1,2))_{1\times(p-1)} & (t_1^2\boldsymbol{\Gamma}(1,2))_{1\times 1} & \cdots & \boldsymbol{0}_{(p-1)\times 1} \\ \vdots & \vdots & \vdots & \vdots & \vdots & \vdots \\ T^2\boldsymbol{\Gamma}(n_z,1)_{(p-1)\times(p-1)} & (Tt_1\boldsymbol{\Gamma}(n_z,1))_{(p-1)\times 1} & T^2\boldsymbol{\Gamma}(n_z,2)_{(p-1)\times(p-1)} & (Tt_1\boldsymbol{\Gamma}(n_z,2))_{(p-1)\times 1} & \cdots & \boldsymbol{0}_{(p-1)\times 1} \\ (Tt_1\boldsymbol{\Gamma}(n_z,1))_{1\times(p-1)} & (t_1^2\boldsymbol{\Gamma}(n_z,1))_{1\times 1} & (Tt_1\boldsymbol{\Gamma}(n_z,2))_{1\times(p-1)} & (t_1^2\boldsymbol{\Gamma}(n_z,2))_{1\times 1} & \cdots & \boldsymbol{0}_{(p-1)\times 1} \\ \boldsymbol{0}_{1\times(p-1)} & \boldsymbol{0}_{1\times(p-1)} & \boldsymbol{0}_{1\times(p-1)} & \boldsymbol{0}_{1\times(p-1)} & \cdots & \boldsymbol{0}_{1\times 1} \end{bmatrix}, \tag{46}$$

$$\begin{cases} \kappa_\alpha \|(\boldsymbol{\Gamma}_3^p)^{\frac{1}{2}} \tilde{\boldsymbol{x}}_3^p\|_2 + \hat{\boldsymbol{d}}_3^p(\tilde{\boldsymbol{x}}_3^p) \leq 0, & \forall (n,p) \in N^2, p > n \\ \kappa_\alpha \|(\boldsymbol{\Gamma}_4^p)^{\frac{1}{2}} \tilde{\boldsymbol{x}}_4^p\|_2 + \hat{\boldsymbol{d}}_4^p(\tilde{\boldsymbol{x}}_4^p) \leq 0, & \forall (n,p) \in N^2, pT > nT + \frac{\chi - \xi}{v_f} \\ \kappa_\alpha \|(\boldsymbol{\Gamma}_5^p)^{\frac{1}{2}} \tilde{\boldsymbol{x}}_5^p\|_2 + \hat{\boldsymbol{d}}_5^p(\tilde{\boldsymbol{x}}_5^p) \leq 0, & \forall (n,p) \in N^2 \\ \qquad \text{s.t. } nT + \frac{\chi - \xi}{v_f} & \in [(p-1)T, pT] \end{cases} \tag{47}$$

$$\begin{aligned} \tilde{\boldsymbol{x}}_3^p := & [\overbrace{q_{\text{out}}(n+1, L_{\text{in}}^z(1)), q_{\text{out}}(n+2, L_{\text{in}}^z(1)), \ldots, q_{\text{out}}(p, L_{\text{in}}^z(1)),}^{p-n} \overbrace{\cdots\cdots}^{(n_z-2)\times(p-n)}, \\ & \overbrace{q_{\text{out}}(n+1, L_{\text{in}}^z(n_z)), q_{\text{out}}(n+2, L_{\text{in}}^z(n_z)), \ldots, q_{\text{out}}(p, L_{\text{in}}^z(n_z))}^{p-n}, 1]^T, \end{aligned} \tag{48}$$

$$\begin{aligned} \hat{\boldsymbol{d}}_3^p := & [\overbrace{T\boldsymbol{P}^z(L_{\text{out}}^z(i), L_{\text{in}}^z(1)), \ldots, T\boldsymbol{P}^z(L_{\text{out}}^z(i), L_{\text{in}}^z(1)),}^{p-n} \overbrace{\cdots\cdots}^{(n_z-2)\times(p-n)}, \\ & \overbrace{T\boldsymbol{P}^z(L_{\text{out}}^z(i), L_{\text{in}}^z(n_z)), \ldots, T\boldsymbol{P}^z(L_{\text{out}}^z(i), L_{\text{in}}^z(n_z))}^{p-n}, -\rho_c v_f[(p-n)T]] \end{aligned} \tag{49}$$

$$\boldsymbol{\Gamma}_3^p := T^2 \begin{bmatrix} \boldsymbol{\Gamma}(1,1)_{(p-n)\times(p-n)} & \boldsymbol{\Gamma}(1,2)_{(p-n)\times(p-n)} & \cdots & \boldsymbol{\Gamma}(1,n_z)_{(p-n)\times(p-n)} & \mathbf{0}_{(p-n)\times 1} \\ \boldsymbol{\Gamma}(2,1)_{(p-n)\times(p-n)} & \boldsymbol{\Gamma}(2,2)_{(p-n)\times(p-n)} & \cdots & \boldsymbol{\Gamma}(2,n_z)_{(p-n)\times(p-n)} & \mathbf{0}_{(p-n)\times 1} \\ \vdots & \vdots & \vdots & \vdots & \vdots \\ \boldsymbol{\Gamma}(n_z,1)_{(p-n)\times(p-n)} & \boldsymbol{\Gamma}(n_z,2)_{(p-n)\times(p-n)} & \cdots & \boldsymbol{\Gamma}(n_z,n_z)_{(p-n)\times(p-n)} & \mathbf{0}_{(p-n)\times 1} \\ \mathbf{0}_{1\times(p-n)} & \mathbf{0}_{1\times(p-n)} & \cdots & \mathbf{0}_{1\times(p-n)} & \mathbf{0}_{1\times 1} \end{bmatrix} \tag{50}$$

$$\begin{aligned} \tilde{\boldsymbol{x}}_4^p := [&\overbrace{q_{\text{out}}(1, L_{\text{in}}^z(1)), q_{\text{out}}(2, L_{\text{in}}^z(1)), \ldots, q_{\text{out}}(n, L_{\text{in}}^z(1)),}^{n} \overbrace{\cdots\cdots}^{(n_z-2)\times n}, \\ &\overbrace{q_{\text{out}}(1, L_{\text{in}}^z(n_z)), q_{\text{out}}(2, L_{\text{in}}^z(n_z)), \ldots, q_{\text{out}}(n, L_{\text{in}}^z(n_z)),}^{n} \\ &\overbrace{q_{\text{out}}(1, L_{\text{out}}^z(i)), q_{\text{out}}(2, L_{\text{out}}^z(i)), \ldots, q_{\text{out}}(p, L_{\text{out}}^z(i)),}^{p} 1]^T, \end{aligned} \tag{51}$$

$$\begin{aligned} \hat{\boldsymbol{d}}_4^p := [&\overbrace{-T\boldsymbol{P}^z(L_{\text{out}}^z(i), L_{\text{in}}^z(1)), \ldots, -T\boldsymbol{P}^z(L_{\text{out}}^z(i), L_{\text{in}}^z(1)),}^{n} \overbrace{\cdots\cdots}^{(n_z-2)\times n}, \\ &\overbrace{-T\boldsymbol{P}^z(L_{\text{out}}^z(i), L_{\text{in}}^z(n_z)), \ldots, -T\boldsymbol{P}^z(L_{\text{out}}^z(i), L_{\text{in}}^z(n_z)),}^{n} \\ &\overbrace{T, \ldots, T}^{p}, -\sum_{k=1}^{k_{max}} \rho(k)X - \rho_c v_f\,(T - t_2)], \end{aligned} \tag{52}$$

$$\boldsymbol{\Gamma}_4^p := T^2 \begin{bmatrix} \boldsymbol{\Gamma}(1,1)_{n\times n} & \boldsymbol{\Gamma}(1,2)_{n\times n} & \cdots & \boldsymbol{\Gamma}(1,n_z)_{n\times n} & \mathbf{0}_{n\times(p+1)} \\ \boldsymbol{\Gamma}(2,1)_{n\times n} & \boldsymbol{\Gamma}(2,2)_{n\times n} & \cdots & \boldsymbol{\Gamma}(2,n_z)_{n\times n} & \mathbf{0}_{n\times(p+1)} \\ \vdots & \vdots & \vdots & \vdots & \vdots \\ \boldsymbol{\Gamma}(n_z,1)_{n\times n} & \boldsymbol{\Gamma}(n_z,2)_{n\times n} & \cdots & \boldsymbol{\Gamma}(n_z,n_z)_{n\times n} & \mathbf{0}_{n\times(p+1)} \\ \mathbf{0}_{(p+1)\times n} & \mathbf{0}_{(p+1)\times n} & \cdots & \mathbf{0}_{(p+1)\times n} & \mathbf{0}_{(p+1)\times(p+1)} \end{bmatrix}, \tag{53}$$

$$\tilde{\boldsymbol{x}}_5^p := \tilde{\boldsymbol{x}}_4^p \tag{54}$$

$$\begin{aligned} \hat{\boldsymbol{d}}_5^p := [&\overbrace{-T\boldsymbol{P}^z(L_{\text{out}}^z(i), L_{\text{in}}^z(1)), \ldots, -T\boldsymbol{P}^z(L_{\text{out}}^z(i), L_{\text{in}}^z(1)),}^{n} \overbrace{\cdots\cdots}^{(n_z-2)\times n}, \\ &\overbrace{-T\boldsymbol{P}^z(L_{\text{out}}^z(i), L_{\text{in}}^z(n_z)), \ldots, -T\boldsymbol{P}^z(L_{\text{out}}^z(i), L_{\text{in}}^z(n_z)),}^{n} \\ &\overbrace{T, \ldots, T}^{p-1}, t_2, -\sum_{k=1}^{k_{max}} \rho(k)X], \end{aligned} \tag{55}$$

$$\boldsymbol{\Gamma}_5^p := \boldsymbol{\Gamma}_4^p, \tag{56}$$

$$\begin{cases} \kappa_\alpha \|(\boldsymbol{\Gamma}_6^p)^{\frac{1}{2}} \tilde{\boldsymbol{x}}_6^p\|_2 + \hat{\boldsymbol{d}}_6^p(\tilde{\boldsymbol{x}}_6^p) \le 0, & \forall (n,p) \in N^2,\ pT > nT + \dfrac{\xi - \chi}{w} \\ \kappa_\alpha \|(\boldsymbol{\Gamma}_7^p)^{\frac{1}{2}} \tilde{\boldsymbol{x}}_7^p\|_2 + \hat{\boldsymbol{d}}_7^p(\tilde{\boldsymbol{x}}_7^p) \le 0, & \forall (n,p) \in N^2 \\ \qquad \text{s.t.}\ \ nT + \dfrac{\xi - \chi}{w} & \in [(p-1)T, pT] \\ M_{\beta_n}(pT, \chi) \ge \beta_p(pT, \chi) & \forall (n,p) \in N^2 \end{cases} \tag{57}$$

$$\begin{aligned} \tilde{\boldsymbol{x}}_6^p := [&\overbrace{q_{\text{out}}(1, L_{\text{in}}^z(1)), q_{\text{out}}(2, L_{\text{in}}^z(1)), \ldots, q_{\text{out}}(n, L_{\text{in}}^z(1)),}^{p} \overbrace{\cdots\cdots}^{(n_z-2)\times p}, \\ &\overbrace{q_{\text{out}}(1, L_{\text{in}}^z(n_z)), q_{\text{out}}(2, L_{\text{in}}^z(n_z)), \ldots, q_{\text{out}}(n, L_{\text{in}}^z(n_z)),}^{p} \\ &\overbrace{q_{\text{out}}(1, L_{\text{out}}^z(i)), q_{\text{out}}(2, L_{\text{out}}^z(i)), \ldots, q_{\text{out}}(p, L_{\text{out}}^z(i)),}^{n} 1]^T, \end{aligned} \tag{58}$$

$$\begin{aligned} \hat{\boldsymbol{d}}_6^p := [&\overbrace{T\boldsymbol{P}^z(L_{\text{out}}^z(i), L_{\text{in}}^z(1)), \ldots, T\boldsymbol{P}^z(L_{\text{out}}^z(i), L_{\text{in}}^z(1)),}^{p} \overbrace{\cdots\cdots}^{(n_z-2)\times p}, \\ &\overbrace{T\boldsymbol{P}^z(L_{\text{out}}^z(i), L_{\text{in}}^z(n_z)), \ldots, T\boldsymbol{P}^z(L_{\text{out}}^z(i), L_{\text{in}}^z(n_z)),}^{p} \\ &\overbrace{-T, \ldots, -T}^{n}, \sum_{k=1}^{k_{max}} \rho(k)X - \rho_c v_f \left((p-n)T - \frac{\xi - \chi}{v_f}\right)] \end{aligned} \tag{59}$$

is a trade-off between the smoothness of the control and the total throughput. The optimal control with $\omega = 20$ is smoother than $\omega = 0.2$ while the overall inflows are lower. Similar to the freeway network in Section IV, Figure 13 shows that the optimal inflows considering uncertainties in turning ratios are lower than the base case, and a higher confidence level induces a larger reduction. It also shows that those optimal inflows are decreasing with time while the control of the base case does not present this trend. However, the optimal inflows of some links, such as links 4 and 5, of the base case in Figure 14 also decrease with time. The reason is that the weight of the total throughput in the objective function, $(n_{max} - i + 1)$, decreases with time step $i$; a large $\omega$ may make the smooth term more significant the throughput term at some point, and the throughput is confined consequently.

## VI. Conclusion

Based on the Lax-Hopf solution to the H-J PDE, this article proposed a robust control model for freeway networks to deal with the uncertainties in turning ratios. The uncertainties are inserted into the model by distributionally robust chance constraints and converted to SOC constraints. Then, multiple case studies for both freeway and urban networks are conducted to investigate the influence of the uncertainties on the interactions between the control of incoming links. To the authors' best knowledge, there are few research studies focusing on the effect of uncertainties in turning ratios on traffic flow control. The proposed model in this article seeks to fit this gap well.

One drawback of this model is that it does not consider the bounds of the turning ratios, which indicates that the real ambiguity set is a subset of the one used in this article. Therefore, this model may provide too conservative of optimal solutions when the variance is large since the worst distribution may not belong to the real ambiguity set. Overcoming this drawback is a promising research direction. Moreover, this article models the uncertainties at different nodes in separate constraints, which implies that the turning ratios at different nodes are independent. How to model them jointly to consider the correlations between intersections is also a meaningful direction. Furthermore, the modeling of uncertainties heavily depends on the type and accuracy of available information about the random parameters. For example, the estimation of the first and the second moments from limited amount of traffic data may have errors. As a result, the actual distribution of the turning ratios may not lie in the used ambiguity set, which can lead to unexpected control consequences. Therefore, it is critical to develop robust control models specific to other types of available information, such as distribution functions, and to investigate the control effect if there exist errors in the prior information. In addition, all model parameters, such as road capacities, are known and fixed in the proposed model. However, in reality and many microscopic traffic simulators such as SUMO and CORSIM, they are usually random and have a significant impact on the control performance. Consequently, the effect of interdependence of different random parameters is another interesting research topic.

## Appendix

### A. Expressions of the SOC Constraints

Equations (40), (47) and (57) are the SOC constraints converted from Equations (15)-(17), respectively, (40), as shown at the page 13, where (41)–(46), as shown at the page 13, where $\mathbf{\Gamma}(i, j)_{a\times b}$ indicates a $a \times b$ matrix in which

$$\mathbf{\Gamma}_6^p := T^2 \begin{bmatrix} \mathbf{\Gamma}(1,1)_{p\times p} & \mathbf{\Gamma}(1,2)_{p\times p} & \cdots & \mathbf{\Gamma}(1,n_z)_{p\times p} & \mathbf{0}_{p\times(n+1)} \\ \mathbf{\Gamma}(2,1)_{p\times p} & \mathbf{\Gamma}(2,2)_{p\times p} & \cdots & \mathbf{\Gamma}(2,n_z)_{p\times p} & \mathbf{0}_{p\times(n+1)} \\ \vdots & \vdots & \vdots & \vdots & \vdots \\ \mathbf{\Gamma}(n_z,1)_{p\times p} & \mathbf{\Gamma}(n_z,2)_{p\times p} & \cdots & \mathbf{\Gamma}(n_z,n_z)_{p\times p} & \mathbf{0}_{p\times(n+1)} \\ \mathbf{0}_{(n+1)\times p} & \mathbf{0}_{(n+1)\times p} & \cdots & \mathbf{0}_{(n+1)\times p} & \mathbf{0}_{(n+1)\times(n+1)} \end{bmatrix}, \tag{60}$$

$$\tilde{\boldsymbol{x}}_7^p := \tilde{\boldsymbol{x}}_6^p, \tag{61}$$

$$\begin{aligned} \hat{\boldsymbol{d}}_7^p := [&\overbrace{T\boldsymbol{P}^z(L_{\text{out}}^z(i), L_{\text{in}}^z(1)), \ldots,}^{p-1} \boldsymbol{P}^z(L_{\text{out}}^z(i), L_{\text{in}}^z(1))t_3, \overbrace{\cdots\cdots}^{(n_z-2)\times p}, \\ &\overbrace{T\boldsymbol{P}^z(L_{\text{out}}^z(i), L_{\text{in}}^z(n_z)), \ldots,}^{p-1} \boldsymbol{P}^z(L_{\text{out}}^z(i), L_{\text{in}}^z(n_z))t_3, \\ &\overbrace{-T, \ldots, -T}^{n}, \sum_{k=1}^{k_{max}} \rho(k)X - \rho_c v_f(\xi - \chi)\left(\frac{1}{w} - \frac{1}{v_f}\right)], \end{aligned} \tag{62}$$

$$\mathbf{\Gamma}_2^p := \begin{bmatrix} T^2\mathbf{\Gamma}(1,1)_{(p-1)\times(p-1)} & (Tt_3\mathbf{\Gamma}(1,1))_{(p-1)\times 1} & T^2\mathbf{\Gamma}(1,2)_{(p-1)\times(p-1)} & (Tt_3\mathbf{\Gamma}(1,2))_{(p-1)\times 1} & \cdots & \mathbf{0}_{(p-1)\times(n+1)} \\ (Tt_3\mathbf{\Gamma}(1,1))_{1\times(p-1)} & (t_3^2\mathbf{\Gamma}(1,1))_{1\times 1} & (Tt_3\mathbf{\Gamma}(1,2))_{1\times(p-1)} & (t_3^2\mathbf{\Gamma}(1,2))_{1\times 1} & \cdots & \mathbf{0}_{(p-1)\times(n+1)} \\ \vdots & \vdots & \vdots & \vdots & \vdots & \vdots \\ T^2\mathbf{\Gamma}(n_z,1)_{(p-1)\times(p-1)} & (Tt_3\mathbf{\Gamma}(n_z,1))_{(p-1)\times 1} & T^2\mathbf{\Gamma}(n_z,2)_{(p-1)\times(p-1)} & (Tt_3\mathbf{\Gamma}(n_z,2))_{(p-1)\times 1} & \cdots & \mathbf{0}_{(p-1)\times(n+1)} \\ (Tt_3\mathbf{\Gamma}(n_z,1))_{1\times(p-1)} & (t_3^2\mathbf{\Gamma}(n_z,1))_{1\times 1} & (Tt_3\mathbf{\Gamma}(n_z,2))_{1\times(p-1)} & (t_3^2\mathbf{\Gamma}(n_z,2))_{1\times 1} & \cdots & \mathbf{0}_{(p-1)\times(n+1)} \\ \mathbf{0}_{(n+1)\times(p-1)} & \mathbf{0}_{(n+1)\times(p-1)} & \mathbf{0}_{(n+1)\times(p-1)} & \mathbf{0}_{(n+1)\times(p-1)} & \cdots & \mathbf{0}_{(n+1)\times(n+1)} \end{bmatrix} \tag{63}$$

all the elements equal $\mathbf{\Gamma}(i, j)$ and $t_1 = \left(\frac{\xi - x_{k-1}}{w} - (p-1)T\right)$, (47), as shown at the page 13, where, (48)–(56), as shown at the pages 13 and 14, where $t_2 = nT + \frac{\chi - \xi}{v_f} - (p-1)T$, (57), as shown at the page 14, where (58)–(63), as shown at the pages 14 and 15.

## References

[1] S. E. Jabari and H. X. Liu, "A stochastic model of traffic flow: Theoretical foundations," *Transp. Res. B, Methodol.*, vol. 46, no. 1, pp. 156–174, Jan. 2012.

[2] X. Qu, J. Zhang, and S. Wang, "On the stochastic fundamental diagram for freeway traffic: Model development, analytical properties, validation, and extensive applications," *Transp. Res. B, Methodol.*, vol. 104, pp. 256–271, Oct. 2017.

[3] J.-B. Sheu and M.-S. Chang, "Stochastic optimal-control approach to automatic incident-responsive coordinated ramp control," *IEEE Trans. Intell. Transp. Syst.*, vol. 8, no. 2, pp. 359–367, Jun. 2007.

[4] B. P. Gokulan and D. Srinivasan, "Distributed geometric fuzzy multiagent urban traffic signal control," *IEEE Trans. Intell. Transp. Syst.*, vol. 11, no. 3, pp. 714–727, Sep. 2010.

[5] B. Heydecker, "Treatment of random variability in traffic modelling," in *Proc. Workshop Traffic Granular Flow HLRZ*. Jülich, Germany: Forschungszentr um Jülich, Oct. 1996, pp. 119–135.

[6] G. F. Newell, "Theory of highway traffic signals," Univ. California, Berkeley, Berkeley, CA, USA, Tech. Rep., 1989.

[7] H. K. Lo, "A reliability framework for traffic signal control," *IEEE Trans. Intell. Transp. Syst.*, vol. 7, no. 2, pp. 250–260, Jun. 2006.

[8] L. Li, W. Huang, and H. K. Lo, "Adaptive coordinated traffic control for stochastic demand," *Transp. Res. C, Emerg. Technol.*, vol. 88, pp. 31–51, Mar. 2018.

[9] M. J. Lighthill and G. B. Whitham, "On kinematic waves. II. A theory of traffic flow on long crowded roads," *Proc. Roy. Soc. London A, Math. Phys. Sci.*, vol. 229, no. 1178, pp. 317–345, 1955.

[10] P. I. Richards, "Shock waves on the highway," *Oper. Res.*, vol. 4, no. 1, pp. 42–51, Feb. 1956.

[11] G. Dervisoglu, G. Gomes, J. Kwon, R. Horowitz, and P. Varaiya, "Automatic calibration of the fundamental diagram and empirical observations on capacity," in *Proc. Transp. Res. Board 88th Annu. Meeting*, vol. 15, 2009, pp. 31–59.

[12] A. Polus and M. A. Pollatschek, "Stochastic nature of freeway capacity and its estimation," *Can. J. Civil Eng.*, vol. 29, no. 6, pp. 842–852, Dec. 2002.

[13] K. Ozbay and E. E. Ozguven, "A comparative methodology for estimating the capacity of a freeway section," in *Proc. IEEE Intell. Transp. Syst. Conf.*, Sep. 2007, pp. 1034–1039.

[14] Y. Li, E. Canepa, and C. Claudel, "Optimal control of scalar conservation laws using linear/quadratic programming: Application to transportation networks," *IEEE Trans. Control Netw. Syst.*, vol. 1, no. 1, pp. 28–39, Mar. 2014.

[15] Y. Li, E. Canepa, and C. Claudel, "Exact solutions to robust control problems involving scalar hyperbolic conservation laws using mixed integer linear programming," in *Proc. 51st Annu. Allerton Conf. Commun., Control, Comput. (Allerton)*, Oct. 2013, pp. 478–485.

[16] H. Liu, C. Claudel, and R. B. Machemehl, "A stochastic formulation of the optimal boundary control problem involving the lighthill whitham richards model," *IFAC-PapersOnLine*, vol. 51, no. 9, pp. 337–342, 2018.

[17] H. Liu, C. Claudel, and R. Machemehl, "Robust traffic control using a first order macroscopic traffic flow model," 2020, *arXiv:2001.06136*. [Online]. Available: http://arxiv.org/abs/2001.06136

[18] G. Como, E. Lovisari, and K. Savla, "Convexity and robustness of dynamic traffic assignment and freeway network control," *Transp. Res. B, Methodol.*, vol. 91, pp. 446–465, Sep. 2016.

[19] J. Ash and D. Newth, "Optimizing complex networks for resilience against cascading failure," *Phys. A, Stat. Mech. Appl.*, vol. 380, pp. 673–683, Jul. 2007.

[20] G. Como, K. Savla, D. Acemoglu, M. A. Dahleh, and E. Frazzoli, "Robust distributed routing in dynamical networks—Part II: Strong resilience, equilibrium selection and cascaded failures," *IEEE Trans. Autom. Control*, vol. 58, no. 2, pp. 333–348, Feb. 2013.

[21] G. Como, K. Savla, D. Acemoglu, M. A. Dahleh, and E. Frazzoli, "Robust distributed routing in dynamical networks—Part I: Locally responsive policies and weak resilience," *IEEE Trans. Autom. Control*, vol. 58, no. 2, pp. 317–332, Feb. 2013.

[22] A. Y. Yazicioglu, M. Roozbehani, and M. A. Dahleh, "Resilient control of transportation networks by using variable speed limits," *IEEE Trans. Control Netw. Syst.*, vol. 5, no. 4, pp. 2011–2022, Dec. 2018.

[23] A. Y. Yazicioglu, M. Roozbehani, and M. A. Dahleh, "Resilient operation of transportation networks via variable speed limits," in *Proc. Amer. Control Conf. (ACC)*, May 2017, pp. 5623–5628.

[24] G. Nilsson, G. Como, and E. Lovisari, "On resilience of multicommodity dynamical flow networks," in *Proc. 53rd IEEE Conf. Decis. Control*, Dec. 2014, pp. 5125–5130.

[25] G. Bianchin, F. Pasqualetti, and S. Kundu, "Resilience of traffic networks with partially controlled routing," in *Proc. Amer. Control Conf. (ACC)*, Jul. 2019, pp. 2670–2675.

[26] R. Arnott, A. de Palma, and R. Lindsey, "Does providing information to drivers reduce traffic congestion?" *Transp. Res. A, Gen.*, vol. 25, no. 5, pp. 309–318, Sep. 1991.

[27] P.-E. Mazaré, A. H. Dehwah, C. G. Claudel, and A. M. Bayen, "Analytical and grid-free solutions to the Lighthill–Whitham–Richards traffic flow model," *Transp. Res. B, Methodol.*, vol. 45, no. 10, pp. 1727–1748, 2011.

[28] V. Henn, "A wave-based resolution scheme for the hydrodynamic LWR traffic flow model," in *Traffic Granular Flow*. Berlin, Germany: Springer, 2005, pp. 105–124.

[29] G. C. Calafiore and L. E. Ghaoui, "On distributionally robust chance-constrained linear programs," *J. Optim. Theory Appl.*, vol. 130, no. 1, pp. 1–22, Dec. 2006.

[30] MOSEK ApS. *The MOSEK Optimization Toolbox for MATLAB Manual, Version 9.0*. Accessed: 2019. [Online]. Available: http://docs.mosek.com/9.0/toolbox/index.html

[31] C. G. Claudel and A. M. Bayen, "Convex formulations of data assimilation problems for a class of Hamilton–Jacobi equations," *SIAM J. Control Optim.*, vol. 49, no. 2, pp. 383–402, Jan. 2011.

[32] E. S. Canepa and C. G. Claudel, "Exact solutions to traffic density estimation problems involving the lighthill-whitham-richards traffic flow model using mixed integer programming," in *Proc. 15th Int. IEEE Conf. Intell. Transp. Syst.*, Sep. 2012, pp. 832–839.

[33] K. Moskowitz, "Discussion of 'freeway level of service as influenced by volume and capacity characteristics' by DR Drew and CJ Keese," *Highway Res. Rec.*, vol. 99, pp. 43–44, Jan. 1965.

[34] E. N. Barron and R. Jensen, "Semicontinuous viscosity solutions for Hamilton–Jacobi equations with convex Hamiltonians," *Commun. Partial Differ. Equ.*, vol. 15, no. 12, pp. 293–309, Jan. 1990.

[35] H. Frankowska, "Lower semicontinuous solutions of Hamilton–Jacobi–Bellman equations," *SIAM J. Control Optim.*, vol. 31, no. 1, pp. 257–272, 1993.

[36] J.-P. Aubin, A. M. Bayen, and P. Saint-Pierre, "Dirichlet problems for some Hamilton–Jacobi equations with inequality constraints," *SIAM J. Control Optim.*, vol. 47, no. 5, pp. 2348–2380, Jan. 2008.

[37] C. G. Claudel and A. M. Bayen, "Lax–Hopf based incorporation of internal boundary conditions into Hamilton–Jacobi equation. Part I: Theory," *IEEE Trans. Autom. Control*, vol. 55, no. 5, pp. 1142–1157, May 2010.

[38] C. G. Claudel and A. M. Bayen, "Lax–Hopf based incorporation of internal boundary conditions into Hamilton-Jacobi equation. Part II: Computational methods," *IEEE Trans. Autom. Control*, vol. 55, no. 5, pp. 1158–1174, May 2010.

[39] R. C. Carlson, I. Papamichail, and M. Papageorgiou, "Local feedback-based mainstream traffic flow control on motorways using variable speed limits," *IEEE Trans. Intell. Transp. Syst.*, vol. 12, no. 4, pp. 1261–1276, Dec. 2011.

[40] M. Gugat, M. Herty, A. Klar, and G. Leugering, "Optimal control for traffic flow networks," *J. Optim. Appl.*, vol. 126, no. 3, pp. 589–616, 2005.

[41] C. Canudas de Wit, "Best-effort highway traffic congestion control via variable speed limits," in *Proc. IEEE Conf. Decis. Control Eur. Control Conf.*, Dec. 2011, pp. 5959–5964.

[42] A. M. Bayen, R. L. Raffard, and C. J. Tomlin, "Network congestion alleviation using adjoint hybrid control: Application to highways," in *Proc. Int. Workshop Hybrid Syst., Comput. Control*. Berlin, Germany: Springer, 2004, pp. 95–110.

[43] E. Delage and Y. Ye, "Distributionally robust optimization under moment uncertainty with application to data-driven problems," *Oper. Res.*, vol. 58, no. 3, pp. 595–612, Jun. 2010.

[44] H. Rahimian and S. Mehrotra, "Distributionally robust optimization: A review," 2019, *arXiv:1908.05659*. [Online]. Available: http://arxiv.org/abs/1908.05659

[45] C. F. Daganzo, "The cell transmission model: A dynamic representation of highway traffic consistent with the hydrodynamic theory," *Transp. Res. B, Methodol.*, vol. 28, no. 4, pp. 269–287, Aug. 1994.

[46] C. F. Daganzo, "The cell transmission model, part II: Network traffic," *Transp. Res. B, Methodol.*, vol. 29, no. 2, pp. 79–93, Apr. 1995.

[47] S. Waller and A. Ziliaskopoulos, "A visual interactive system for transportation algorithms," in *Proc. 78th Annu. Meeting Transp. Res. Board*, Washington, DC, USA, 1998.

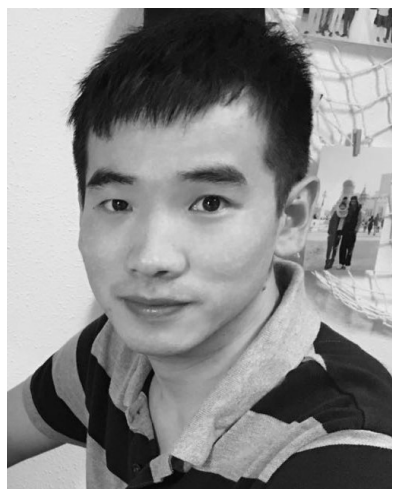

**Hao Liu** received the M.S. degree in solid mechanics from Peking University in 2015, the M.S. degree in statistics from The University of Texas at Austin, USA, in 2018, and the Ph.D. degree in civil, architectural and environmental engineering from The University of Texas at Austin. His research interests include traffic flow modeling and optimization.

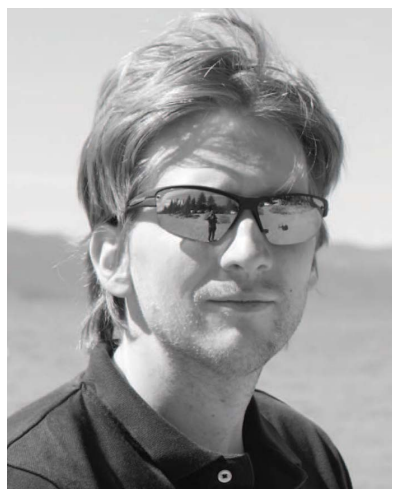

**Christian Claudel** received the M.S. degree in plasma physics from the École Normale Supérieure de Lyon in 2004 and the Ph.D. degree in EECS from UC Berkeley in 2010. He was an Assistant Professor of electrical engineering with the King Abdullah University of Science and Technology. He is currently an Assistant Professor of civil architectural and environmental engineering with The University of Texas at Austin. His research interests include the control and estimation of distributed parameter systems, wireless sensor networks, and unmanned air vehicles.

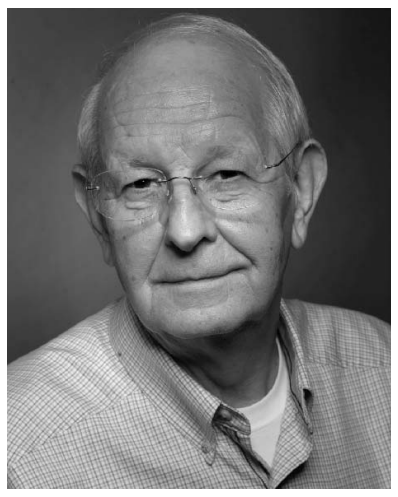

**Randy Machemehl** is currently a Professor of transportation engineering with The University of Texas at Austin and the former Director of the Center for Transportation Research. He was in private engineering practice as a Staff Member of Wilbur Smith and an Associate before joining the faculty of The University of Texas at Austin in 1978. His research interests include transportation system operations, public transportation systems planning and design, traffic data acquisition, traffic simulation, optimization, and bicycle safety.

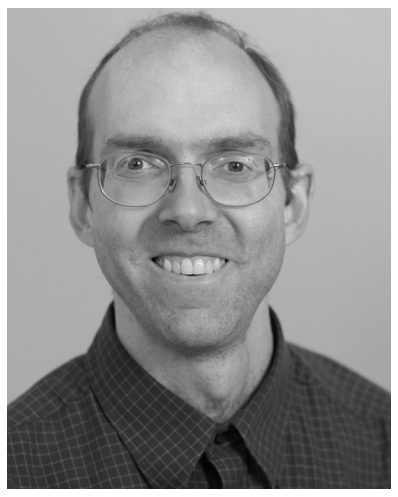

**Kenneth A. Perrine** received the M.S. degree in civil engineering from the University of Washington. He has worked with the Pacific Northwest National Laboratory on high-performance computing. He is currently a Research Associate with the Network Modeling Center, Center for Transportation Research, The University of Texas at Austin. His research interests include traffic modeling, data analysis, and sensors.